\documentclass[reqno,12pt,a4paper]{amsart}

\usepackage{cancel}
\usepackage{mathrsfs}
\usepackage{amssymb}
\usepackage{amsfonts}
\usepackage{latexsym}
\usepackage{amsthm}
\usepackage[numbers, sort&compress]{natbib}
\usepackage{xcolor}
\usepackage{graphicx}
\usepackage[T2A,T1]{fontenc}
\usepackage[utf8]{inputenc}
\usepackage[russian,english]{babel}

\def\lb{\label}

\newcommand{\er}[1]{\textrm{(\ref{#1})}}

\begin{document}


\renewcommand{\theequation}{\arabic{section}.\arabic{equation}}
\theoremstyle{plain}
\newtheorem{theorem}{\bf Theorem}[section]
\newtheorem{lemma}[theorem]{\bf Lemma}
\newtheorem{corollary}[theorem]{\bf Corollary}
\newtheorem{proposition}[theorem]{\bf Proposition}
\newtheorem{definition}[theorem]{\bf Definition}

\def\a{\alpha}  \def\cA{{\mathcal A}}     \def\bA{{\bf A}}  \def\mA{{\mathscr A}}
\def\b{\beta}   \def\cB{{\mathcal B}}     \def\bB{{\bf B}}  \def\mB{{\mathscr B}}
\def\g{\gamma}  \def\cC{{\mathcal C}}     \def\bC{{\bf C}}  \def\mC{{\mathscr C}}
\def\G{\Gamma}  \def\cD{{\mathcal D}}     \def\bD{{\bf D}}  \def\mD{{\mathscr D}}
\def\d{\delta}  \def\cE{{\mathcal E}}     \def\bE{{\bf E}}  \def\mE{{\mathscr E}}
\def\D{\Delta}  \def\cF{{\mathcal F}}     \def\bF{{\bf F}}  \def\mF{{\mathscr F}}
\def\c{\chi}    \def\cG{{\mathcal G}}     \def\bG{{\bf G}}  \def\mG{{\mathscr G}}
\def\z{\zeta}   \def\cH{{\mathcal H}}     \def\bH{{\bf H}}  \def\mH{{\mathscr H}}
\def\e{\eta}    \def\cI{{\mathcal I}}     \def\bI{{\bf I}}  \def\mI{{\mathscr I}}
\def\p{\psi}    \def\cJ{{\mathcal J}}     \def\bJ{{\bf J}}  \def\mJ{{\mathscr J}}
\def\vT{\Theta} \def\cK{{\mathcal K}}     \def\bK{{\bf K}}  \def\mK{{\mathscr K}}
\def\k{\kappa}  \def\cL{{\mathcal L}}     \def\bL{{\bf L}}  \def\mL{{\mathscr L}}
\def\l{\lambda} \def\cM{{\mathcal M}}     \def\bM{{\bf M}}  \def\mM{{\mathscr M}}
\def\L{\Lambda} \def\cN{{\mathcal N}}     \def\bN{{\bf N}}  \def\mN{{\mathscr N}}
\def\m{\mu}     \def\cO{{\mathcal O}}     \def\bO{{\bf O}}  \def\mO{{\mathscr O}}
\def\n{\nu}     \def\cP{{\mathcal P}}     \def\bP{{\bf P}}  \def\mP{{\mathscr P}}
\def\r{\rho}    \def\cQ{{\mathcal Q}}     \def\bQ{{\bf Q}}  \def\mQ{{\mathscr Q}}
\def\s{\sigma}  \def\cR{{\mathcal R}}     \def\bR{{\bf R}}  \def\mR{{\mathscr R}}
\def\vs{\varsigma}  \def\cS{{\mathcal S}}     \def\bS{{\bf S}}  \def\mS{{\mathscr S}}
\def\t{\tau}    \def\cT{{\mathcal T}}     \def\bT{{\bf T}}  \def\mT{{\mathscr T}}
\def\f{\phi}    \def\cU{{\mathcal U}}     \def\bU{{\bf U}}  \def\mU{{\mathscr U}}
\def\F{\Phi}    \def\cV{{\mathcal V}}     \def\bV{{\bf V}}  \def\mV{{\mathscr V}}
\def\P{\Psi}    \def\cW{{\mathcal W}}     \def\bW{{\bf W}}  \def\mW{{\mathscr W}}
\def\o{\omega}  \def\cX{{\mathcal X}}     \def\bX{{\bf X}}  \def\mX{{\mathscr X}}
\def\x{\xi}     \def\cY{{\mathcal Y}}     \def\bY{{\bf Y}}  \def\mY{{\mathscr Y}}
\def\X{\Xi}     \def\cZ{{\mathcal Z}}     \def\bZ{{\bf Z}}  \def\mZ{{\mathscr Z}}
\def\O{\Omega}


\def\mb{{\mathscr b}}
\def\mh{{\mathscr h}}
\def\me{{\mathscr e}}
\def\mk{{\mathscr k}}
\def\mz{{\mathscr z}}
\def\mx{{\mathscr x}}

\def\be{{\bf e}} \def\bc{{\bf c}} \def\bt{{\bf t}}
\def\bx{{\bf x}} \def\by{{\bf y}}
\def\bv{{\bf v}} 
\def\Om{\Omega} \def\bp{{\bf p}}

\newcommand{\gA}{\mathfrak{A}}          \newcommand{\ga}{\mathfrak{a}}
\newcommand{\gB}{\mathfrak{B}}          \newcommand{\gb}{\mathfrak{b}}
\newcommand{\gC}{\mathfrak{C}}          \newcommand{\gc}{\mathfrak{c}}
\newcommand{\gD}{\mathfrak{D}}          \newcommand{\gd}{\mathfrak{d}}
\newcommand{\gE}{\mathfrak{E}}
\newcommand{\gF}{\mathfrak{F}}           \newcommand{\gf}{\mathfrak{f}}
\newcommand{\gG}{\mathfrak{G}}           \newcommand{\gog}{\mathfrak{g}}
\newcommand{\gH}{\mathfrak{H}}           \newcommand{\gh}{\mathfrak{h}}
\newcommand{\gI}{\mathfrak{I}}           \newcommand{\gi}{\mathfrak{i}}
\newcommand{\gJ}{\mathfrak{J}}           \newcommand{\gj}{\mathfrak{j}}
\newcommand{\gK}{\mathfrak{K}}            \newcommand{\gk}{\mathfrak{k}}
\newcommand{\gL}{\mathfrak{L}}            \newcommand{\gl}{\mathfrak{l}}
\newcommand{\gM}{\mathfrak{M}}            \newcommand{\gm}{\mathfrak{m}}
\newcommand{\gN}{\mathfrak{N}}            \newcommand{\gn}{\mathfrak{n}}
\newcommand{\gO}{\mathfrak{O}}
\newcommand{\gP}{\mathfrak{P}}             \newcommand{\gp}{\mathfrak{p}}
\newcommand{\gQ}{\mathfrak{Q}}             \newcommand{\gq}{\mathfrak{q}}
\newcommand{\gR}{\mathfrak{R}}             \newcommand{\gr}{\mathfrak{r}}
\newcommand{\gS}{\mathfrak{S}}              \newcommand{\gs}{\mathfrak{s}}
\newcommand{\gT}{\mathfrak{T}}             \newcommand{\gt}{\mathfrak{t}}
\newcommand{\gU}{\mathfrak{U}}             \newcommand{\gu}{\mathfrak{u}}
\newcommand{\gV}{\mathfrak{V}}             \newcommand{\gv}{\mathfrak{v}}
\newcommand{\gW}{\mathfrak{W}}             \newcommand{\gw}{\mathfrak{w}}
\newcommand{\gX}{\mathfrak{X}}               \newcommand{\gx}{\mathfrak{x}}
\newcommand{\gY}{\mathfrak{Y}}              \newcommand{\gy}{\mathfrak{y}}
\newcommand{\gZ}{\mathfrak{Z}}             \newcommand{\gz}{\mathfrak{z}}

\def\ve{\varepsilon} \def\vt{\vartheta} \def\vp{\varphi}  \def\vk{\varkappa}
\def\vr{\varrho}

\def\A{{\mathbb A}} \def\B{{\mathbb B}} \def\C{{\mathbb C}}
\def\dD{{\mathbb D}} \def\E{{\mathbb E}} \def\dF{{\mathbb F}} \def\dG{{\mathbb G}}
\def\H{{\mathbb H}}\def\I{{\mathbb I}} \def\J{{\mathbb J}} \def\K{{\mathbb K}}
\def\dL{{\mathbb L}}\def\M{{\mathbb M}} \def\N{{\mathbb N}} \def\dO{{\mathbb O}}
\def\dP{{\mathbb P}} \def\dQ{{\mathbb Q}} \def\R{{\mathbb R}}\def\S{{\mathbb S}} \def\T{{\mathbb T}}
\def\U{{\mathbb U}} \def\V{{\mathbb V}}\def\W{{\mathbb W}} \def\X{{\mathbb X}}
\def\Y{{\mathbb Y}} \def\Z{{\mathbb Z}}

\def\dk{{\Bbbk}}


\def\la{\leftarrow}              \def\ra{\rightarrow}     \def\Ra{\Rightarrow}
\def\ua{\uparrow}                \def\da{\downarrow}
\def\lra{\leftrightarrow}        \def\Lra{\Leftrightarrow}
\newcommand{\abs}[1]{\lvert#1\rvert}
\newcommand{\br}[1]{\left(#1\right)}

\def\lan{\langle} \def\ran{\rangle}


\def\lt{\biggl}                  \def\rt{\biggr}
\def\ol{\overline}               \def\wt{\widetilde}
\def\no{\noindent}


\let\ge\geqslant                 \let\le\leqslant
\def\lan{\langle}                \def\ran{\rangle}
\def\/{\over}                    \def\iy{\infty}
\def\sm{\setminus}               \def\es{\emptyset}
\def\ss{\subset}                 \def\ts{\times}
\def\pa{\partial}                \def\os{\oplus}
\def\om{\ominus}                 \def\ev{\equiv}
\def\iint{\int\!\!\!\int}        \def\iintt{\mathop{\int\!\!\int\!\!\dots\!\!\int}\limits}
\def\wh{\widehat}
\def\bs{\backslash}
\def\na{\nabla}
\def\ti{\tilde}
\def\hb{\hbar}
\def\cd{\centerdot}
\def\we{\wedge}

\def\sh{\mathop{\mathrm{sh}}\nolimits} \def\ch{\mathop{\mathrm{ch}}\nolimits}
\def\all{\mathop{\mathrm{all}}\nolimits}
\def\Area{\mathop{\mathrm{Area}}\nolimits}
\def\arg{\mathop{\mathrm{arg}}\nolimits}
\def\const{\mathop{\mathrm{const}}\nolimits}
\def\det{\mathop{\mathrm{det}}\nolimits}
\def\diag{\mathop{\mathrm{diag}}\nolimits}
\def\diam{\mathop{\mathrm{diam}}\nolimits}
\def\dim{\mathop{\mathrm{dim}}\nolimits}
\def\dist{\mathop{\mathrm{dist}}\nolimits}
\def\Im{\mathop{\mathrm{Im}}\nolimits}
\def\Iso{\mathop{\mathrm{Iso}}\nolimits}
\def\Ker{\mathop{\mathrm{Ker}}\nolimits}
\def\Lip{\mathop{\mathrm{Lip}}\nolimits}
\def\rank{\mathop{\mathrm{rank}}\limits}
\def\Ran{\mathop{\mathrm{Ran}}\nolimits}
\def\Re{\mathop{\mathrm{Re}}\nolimits}
\def\Res{\mathop{\mathrm{Res}}\nolimits}
\def\res{\mathop{\mathrm{res}}\limits}
\def\sign{\mathop{\mathrm{sign}}\nolimits}
\def\span{\mathop{\mathrm{span}}\nolimits}
\def\supp{\mathop{\mathrm{supp}}\nolimits}
\def\Tr{\mathop{\mathrm{Tr}}\nolimits}
\def\BBox{\hspace{1mm}\vrule height6pt width5.5pt depth0pt \hspace{6pt}}
\def\where{\mathop{\mathrm{where}}\nolimits}
\def\as{\mathop{\mathrm{as}}\nolimits}
\def\1{1\!\!1}



\newcommand\nh[2]{\widehat{#1}\vphantom{#1}^{(#2)}}
\def\di{\diamond}

\def\Oplus{\bigoplus\nolimits}



\def\qqq{\qquad}
\def\qq{\quad}
\let\ge\geqslant
\let\le\leqslant
\let\geq\geqslant
\let\leq\leqslant
\newcommand{\ca}{\begin{cases}}
\newcommand{\ac}{\end{cases}}
\newcommand{\ma}{\begin{pmatrix}}
\newcommand{\am}{\end{pmatrix}}
\renewcommand{\[ }{\begin{equation}}
\renewcommand{\] }{\end{equation}}
\def\eq{\begin{equation}}
\def\qe{\end{equation}}
\def\[{\begin{equation}}
\def\bu{\bullet}
\def\ced{\centerdot}
\def\tes{\textstyle}
\def\ci{\circ}


\title[{Inverse problems for ZS operators and their isomorphisms, II
 }]
{Inverse problems for Zakharov--Shabat operators and their
isomorphisms, II}

\date{\today}

\author[Evgeny Korotyaev]{Evgeny Korotyaev}
\author[Zongfeng Zhang]{Zongfeng Zhang}

\address{Department of Analysis, Saint-Petersburg State University,
 Universitetskaya nab. 7/9, St. Petersburg, 199034, Russia, \ korotyaev@gmail.com}

\address{Academy for Advanced Interdisciplinary Studies, Northeast Normal University,
 Changchun, China, \ zhangzf@nenu.edu.cn}


\begin{abstract}
We study isomorphic inverse spectral problems for Zakharov--Shabat
operators on the unit interval with generic boundary conditions. A
unitary gauge transformation reduces these problems to the
corresponding problems with basic boundary conditions, with an
explicit shift of the spectral parameter. Using previous results for
canonical boundary conditions and four-spectra data  mapping we
describe isomorphic inverse spectral problems for Zakharov--Shabat
operators. Moreover, we study other inverse problems including
replacing mappings and canonical relations.
\end{abstract}

\maketitle

\begin{quotation}
\begin{center}
{\bf Table of Contents}
\end{center}

\vskip 6pt

{\footnotesize

1. Introduction and main results \hfill \pageref{Sec1}\ \ \ \ \

2. Preliminary results \hfill\pageref{Sec2}\ \ \ \ \

3. Gauge transformations and their properties  \hfill \pageref{Sec3}\ \ \ \ \

4. Isomorphic inverse problems on the unit interval \hfill \pageref{Sec4}\ \ \ \ \

 }
\end{quotation}

\section {Introduction and main results \lb{Sec1}}
\setcounter{equation}{0}
\subsection{Introduction}
We study isomorphic inverse spectral problems for Zakharov--Shabat
systems (ZS systems) $Jy'+Vy=\l y$ on the unit interval with generic
boundary conditions.  Inverse spectral problems have been
extensively studied, with a rich collection of approaches and
results available in \cite{AT04,FT86,L87,LS91,M86,PT87} and the
references therein. In general, the study of an inverse spectral
problem consists of the following parts:

\no (i) Uniqueness: prove that the spectral data (eigenvalues together with
certain additional parameters) determine the potential uniquely. \\
(ii) Reconstruction: reconstruct the potential from spectral data. \\
(iii) Characterization: describe all spectral data corresponding to fixed classes of potentials. \\
(iv) Stability: obtain a priori two--sided estimates relating the potential and spectral data. \\
(v) Isomorphism: prove that the spectral mappings are isomorphic.

Our main goal is to discuss the last property, namely the {\it
isomorphism} of inverse problems. Isomorphisms between inverse
spectral problems have been studied for Sturm--Liouville operators
on the unit interval and on the circle \cite{K25}, for Jacobi
operators \cite{KL27}, and for ZS operators with basic boundary
conditions \cite{KZ26}. Related questions for operators in impedance
form are considered in the work in preparation \cite{K27}. The
present paper extends the ZS results to separated generic boundary
conditions. Using previous results for canonical boundary conditions
and four-spectra data  mapping we describe isomorphic inverse
spectral problems for Zakharov--Shabat operators. The key step is a
gauge transformation that rotates the potential and shifts every
eigenvalue by the same constant. This reduces the analytic
parametrization results to those in \cite{KZ26} and yields explicit
transformations between the resulting inverse problems.

\begin{definition}
Let $f,f_0:\cK\to X$ be mappings from a real Hilbert space $\cK$ to
a set $X$. They are called \emph{isomorphic} if $f_0=f\ci U$ for
some (not necessarily) isomorphism $U$ of $\cK$ onto itself.
\end{definition}

Note that if $f$ is bijective, the intertwining map is uniquely
determined by $U=f^{-1}\circ f_0$. Isomorphic mappings have the same
image and are injective or bijective simultaneously. Analytic
parametrization results also transfer when the target has the
corresponding analytic structure. An identity for $f$ transfers to
$f_0$ after replacing the potential $v$ by $Uv$. In particular, a
priori bounds in the original potential norm transfer unchanged when
$U$ preserves that norm. More general stability estimates require
quantitative control of $U$ and $U^{-1}$.

We consider the inverse spectral problems for the ZS system
\[
\lb{eqZS}
Jy'+Vy=\l y
\]
on the interval $[0,1]$, subject to four distinct generic boundary conditions:
\[
\lb{sgc1}
\tes y_1(0)=0, \qqq y_1(1)c_\a  -y_2(1)s_\a  =0,
\]
\[
\lb{sgc2}
\tes y_2(0)=0, \qqq y_1(1)s_\a  +y_2(1)c_\a  =0,
\]
\[
\lb{sgc3}
\tes y_1(0)=0, \qqq y_1(1)s_\a  +y_2(1)c_\a  =0,
\]
\[
\lb{sgc4}
\tes y_2(0)=0, \qqq y_1(1)c_\a  -y_2(1)s_\a  =0.
\]
The matrix-valued potential $V$ and the matrix $J$ are given by
$$
\begin{aligned}
V=\ma v_1&v_2\\v_2&-v_1\am, \qqq  J=\ma0&1\\-1&0\am.
\end{aligned}
$$
Here $v=(v_1,v_2)^\top$ belongs to the real Hilbert space
$\mH=L^2((0,1),\R)\os L^2((0,1),\R)$ equipped with norm
$\|v\|^2=\int_0^1(v_1^2+v_2^2)\,dx$. Throughout, we use the notation
$c_\a =\cos\a$, $s_\a =\sin\a$ and assume that
$$
\tes \a\in\A:=\big(-\frac{\pi}{2},\frac{\pi}{2}\big).
$$
For $\a=0$, problems \er{sgc1}--\er{sgc4} reduce successively to the
Dirichlet problem, the Neumann problem, and two mixed problems. For
$\a=\pm\frac{\pi}{2}$, the same four types of boundary conditions
occur, but in different order. These special cases, $\a=0,
\pm\frac{\pi}{2}$, are discussed in \cite{KZ26}. In the present
paper, we consider the general case $\a\in\A$.

\no {\bf Remark on boundary conditions.} A general separated
boundary value problem on $[0,1]$  consists of the ZS system
together with the boundary conditions
$$
Jy'+Vy=\l y, \qqq \left\{\begin{aligned} y_1(0)c_\b-y_2(0)s_\b=0
\\
y_1(1)c_\a -y_2(1)s_\a =0\end{aligned}\right. , \qqq \a, \b\in\A.
$$
Lemma \ref{GTsG} reduces this problem to \er{sgc1} under the gauge
with a rotated potential and spectral parameter $\l-\b$. Thus the
operator is unitarily equivalent to the reduced operator plus
$\b\1_2$; the spectral shift is part of the equivalence. Analogous
reductions apply to the other three basic choices at the left
endpoint.

Let $Y(x,\l,v)=\big(\vt,\vp\big)(x,\l,v)$ denote the fundamental
matrix  solution of the system \er{eqZS}, {{normalized by
$Y(0,\l,v)=\1_2$}}, where $\1_2$ is the identity $2\ts 2$ matrix. Its
columns are the vector-valued functions $\vt=(\vt_1, \vt_2)^\top$
and $\vp=(\vp_1, \vp_2)^\top$, which satisfy the initial conditions
$$
\vt(0,\l,v)=(1, 0)^\top, \qqq \vp(0,\l,v)=(0, 1)^\top.
$$
For $n\in\Z$, let $\m_{\a,n}$, $\n_{\a,n}$, $\t_{\a,n}$, and
$\vr_{\a,n}$ denote the corresponding eigenvalues of the problems
\er{sgc1}--\er{sgc4}, respectively. These eigenvalues are known to
be simple and to satisfy the following interlacing relations:
$$
\begin{aligned}
& ...<\ol {\t_{\a,1}, \vr_{\a,1}}<\ol {\mu_{\a,1},\nu_{\a,1}}
<\ol {\t_{\a,2}, \vr_{\a,2}}< \ol
{\mu_{\a,2},\nu_{\a,2}}<...,
\end{aligned}
$$
where, for brevity, $\ol {u,v}$ denotes $\min \{u,v\}\le \max \{u,v\}$.
The eigenvalues have the following asymptotics:
$$
\begin{aligned}
& \n_{\a,n},\m_{\a,n}=\mu_{\a,n}^o+o(1),\qq \t_{\a,n},
\vr_{\a,n}=\tau_{\a,n}^o+o(1) \qq \as \ n\to \pm \infty,
\end{aligned}
$$
where the corresponding unperturbed eigenvalues satisfy
$$
\tes \nu_{\a,n}^o=\mu_{\a,n}^o=\pi n-\a, \qqq
\tau_{\a,n}^o=\vr_{\a,n}^o = \pi\big(n-\frac{1}{2}\big)-\a, \qqq
n\in \Z.
$$

\subsection{Main results}
Let  $\cK$ and $\cS$ be real separable Hilbert spaces. Denote  by
$\cK_\C$ and $\cS_\C$ their respective complexifications. A mapping
$f: \cK\to\cS$ is a local real analytic isomorphism if and only if,
for every $y\in\cK$, $f$ admits an analytic continuation $\wt f$
into a complex neighborhood $\V\ss\cK_\C$ of $y$ such that $\wt f$
is a bijection between $\V$ and some open set $\wt f(\V)\ss\cS_\C$,
and both $\wt f$ and $\wt f^{-1}$ are analytic on $\V$ and $\wt
f(\V)$, respectively. A mapping $f$ is a real-analytic bijection
(RAB) between $\cK$ and $\cS$ if it is both a bijection and a local
real analytic isomorphism.

For $1\le p<\infty$, let $\ell^p=\ell^p(\Z)$ be the real Banach
space of sequences $f=(f_n)_{n\in\Z}$ equipped with the norm
$\|f\|_{(p)}^p=\sum_{n\in\Z}|f_n|^p<\infty$. We write
$\|f\|=\|f\|_{(2)}$ when $p=2$. To formulate our main results, we
combine the four one-spectra mappings  into a map of spectral
differences. Following earlier work \cite{K25, K19}, the
four-spectra mapping for ZS operators with $\a=0$ was introduced in
\cite{KZ26}. We define the four-spectra mapping $\gf_\a : \mH\to
\ell^2$ by
\begin{equation}
\label{Def4sma} \tes
 v\to \gf_\a (v)=(\gf_{\a,n}(v))_{n\in \Z},\qq
\gf_{\a,2n-1}=\frac{1}{2}(\vr_{\a,n} -\t_{\a,n} ),\qqq
\gf_{\a,2n}=\frac{1}{2}(\n_{\a,n} -\m_{\a,n}), \ \ \ \ \ n\in \Z.
\end{equation}
We show that the four-spectra mapping is a real-analytic bijection
between  the space of potentials and the corresponding spectral
data.

\begin{theorem}
\label{T4sma} Let $\a\in\A$. The four-spectra mapping $\gf_\a :
\mH\to \ell^2$ defined by \er{Def4sma}, is a RAB from $\mH$ onto
$\ell^2$ and  satisfies
\begin{equation}
\label{egf1}
\begin{aligned}
\tes \frac{1}{\sqrt 2}\|\gf_\a (v)\|\leq \|v\|\leq 2\|\gf_\a (v)\|
\big(1+\|\gf_\a (v)\|\big), \qqq \forall \ v\in\mH.
\end{aligned}
\end{equation}
\end{theorem}

\no {\bf Remark.} For fixed $\a$, let $\mD_\a$ consist of those
$v\in\mH$ for which $\gf_\a(v)$ is finitely supported. Since
finitely supported sequences are dense in $\ell^2$ and $\gf_\a^{-1}$
is continuous, $\mD_\a$ is dense in $\mH$. The notation
distinguishes this set from the basic set $\mD$ used in Section~2.

Following P\"oschel and Trubowitz \cite{PT87}, we introduce the
following sets of strictly increasing real sequences:
$$
\begin{aligned}
&\gJ^o_\a=\Big\{ s\!=\!(s_n)_{n\in \Z}: ....<s_1<s_2<....., \qq s_n\!
=\!\pi n-\a\!+\breve{s}_{n\,},
 \qqq \breve{s}=(\breve{s}_{n\,})_{n\in \Z}\!\in\!\ell^2 \Big\},
 \\
&\tes \gJ^1_\a=\Big\{ s\!=\!(s_n)_{n\in \Z}: ....<s_1<s_2<....., \qq
s_n\!=\!\pi (n-\frac{1}{2})-\a\!+\breve{s}_{n\,}, \qqq \breve{s}
=(\breve{s}_{n\,})_{n\in \Z}\!\in\!\ell^2 \Big\},
\\
&  \tes \gJ_\a=\Big\{ t\!=\!(t_n)_{n\in \Z}: ...<t_1<t_2<....., \qq
t_n\!=\frac{\pi n}{2}-\a+\breve{t}_{n\,}, \qqq \breve{t}=
(\breve{t}_{n\,})_{n\in \Z}\!\in\!\ell^2 \Big\}.
\end{aligned}
$$
The mapping $s=(s_n)_{n\in \Z}\lra \breve{s} $ is a natural
coordinate map between $\gJ^o_\a$ and some open convex subset
$\breve{\gJ}^o_\a=\Big\{ \breve{s}\!=\!(\breve{ s}_n)_{n\in
\Z}\in\!\ell^2:
...<\mu_{\a,1}^o+\breve{s}_{1\,}<\mu_{\a,2}^o+\breve{s}_{2\,}<.....
\Big\} $ of $\ell^2$\,. Following  \cite{PT87} we identify
$\gJ^o_\a$ and $\breve{\gJ}^o_\a$ using this mapping.  Below we
refer to $\breve{s}=\!(\breve{ s}_n)_{n\in \Z}\in\!\ell^2$ as the
standard coordinate system on $\gJ^o_\a$. As in \cite{PT87} this
identification allows us to carry out analysis on $\gJ^o_\a$ as if
it were an open convex subset of $\ell^2$. Similar standard
coordinate systems can be defined on $\gJ^1_\a$ and $\gJ_\a$. Note
that for $\a=0$, $\gJ^o_\a$, $\gJ^1_\a$, $\gJ_\a$ reduce to $\gJ^o$,
$\gJ^1$, $\gJ$, respectively.

To describe the inverse problems, we introduce one-spectra mappings
$\m_\a, \nu_\a: \mH\to \gJ^o_\a$, and $\t_\a, \vr_\a: \mH\to
\gJ^1_\a$ by
$$
\label{mntz1}
\begin{aligned}
v\to \m_\a=(\m_{\a,n})_{n\in \Z}, \qq  v\to \nu_\a=(\n_{\a,n})_{n\in
\Z}, \qq v\to \t_\a=(\t_{\a,n})_{n\in \Z}, \qq v\to
\vr_\a=(\vr_{\a,n})_{n\in \Z}.
\end{aligned}
$$
Given two strictly increasing, interlacing sequences of eigenvalues,
we define the associated two-spectra mapping. For example, for
$\t_\a=(\tau_{\a,n})_{n\in \Z}\in \gJ^1_\a$ and
$\n_\a=(\n_{\a,n})_{n\in \Z} \in \gJ^o_\a$ satisfying
$$
...<\tau_{\a,1}<\n_{\a,1}<\tau_{\a,2}<\n_{\a,2}<...,
$$
the two-spectra mapping $\t_\a \star \n_\a: \mH \to \gJ_\a$ is given by
$$
(\t_\a\star\n_\a)_{2n-1}=\t_{\a,n},\qquad
(\t_\a\star\n_\a)_{2n}=\n_{\a,n},\qquad n\in\Z.
$$
The other interlaced mappings use the same odd/even indexing. Below
we need multiplication operators $F_\a$ and $\1_-$ on $\mH$ given by
\[
\lb{DefF_a} F_\a v=e^{\a xJ} v,\qqq \1_-v=-v,\qqq \a\in \A.
\]
 We define the reflection operator $\cR$ and the reflection-type
operators $\mU_1, \mU_2$ acting on $\mH$ by
\[
\lb{DefF_12} (\cR v)(x)= v(1-x), \qq x\in[0,1], \qqq \mU_1=\ma
c_{2\a} & -s_{2\a}  \\
-s_{2\a} & -c_{2\a} \am\cR, \qqq
\mU_2=-\mU_1.
\]
Finally, define the shifted operator $\cS_\a: \gJ\to \gJ_\a$ by
\[
\lb{DefShiftO}
(\cS_\a a)_n=a_n-\a,\qquad
\cS_\a(a,b)=(\cS_\a a,\cS_\a b).
\]
We use the same notation  for the operators $\gJ^o\to\gJ^o_\a$ and
$\gJ^1\to\gJ^1_\a$, and apply it componentwise to tuples of spectra.
In standard sequence coordinates, each operator is the identity, so
it is a real-analytic isomorphism with inverse $\cS_{-\a}$.

We establish isomorphisms between two-spectra mappings.

\begin{theorem} \lb{Tiso2sm} Let $\a\in\A$.
Then all two-spectra mappings $\t_\a\star \m_\a$, $\t_\a\star
\n_\a$, $\vr_\a\star \m_\a$, $\vr_\a\star \n_\a $ acting  from $\mH$
to $\gJ_\a$ are pairwise isomorphic. Each of them is a RAB from
$\mH$ onto $\gJ_\a$ and satisfies
\[
\lb{2sma}
\begin{aligned}
\t_\a\star \m_\a
=(\vr_\a\star \n_\a )\circ\1_-=(\vr_\a\star \m_\a)\circ\mU_1
=(\t_\a\star \n_\a )\circ\mU_2=\cS_\a(\t\star \m)\ci F_{2\a}^*.
\end{aligned}
\]
\end{theorem}

\no Note that  the authors proved these results for the case $\a=0$
in \cite{KZ26}.

For each eigenvalue $\m_{\a,n}$,  $\vp(x,\m_{\a,n},v)$ is the
associated eigenfunction. It is known that the vector
$\vp(1,\m_{\a,n},v)$ satisfies
$$
(-1)^ne^{\gr_{\a,n}(v)}\vp(1,\m_{\a,n},v)=(\sin\a,\cos\a)^\top.
$$
In view of these normalization conditions for the respective
boundary  value problems, we define the norming constants
$\gr_{\a,n}$, $\gs_{\a,n}$, $\gt_{\a,n}$, $\gu_{\a,n}$ and the
associated norming mappings by
\[
\lb{DefNorming}
\begin{aligned}
&\gr_{\a,n}=-\ln\Big|\frac{\vp_2(1,\m_{\a,n},v)}{c_\a}\Big|, \qqq
\gs_{\a,n}=-\ln\Big|\frac{\vt_1(1,\n_{\a,n},v)}{c_\a}\Big|,  \\
&\gt_{\a,n}=-\ln\Big|\frac{\vp_1(1,\t_{\a,n},v)}{c_\a}\Big|, \qqq
\gu_{\a,n}=-\ln\Big|\frac{\vt_2(1,\vr_{\a,n},v)}{c_\a}\Big|,
\end{aligned}
\]
and
$$
\begin{aligned}
v\to \gr_\a =(\gr_{\a,n} )_{n\in \Z},\qq
v\to \gs_\a =(\gs_{\a,n} )_{n\in \Z}, \qq
v\to \gt_\a =(\gt_{\a,n} )_{n\in \Z},\qq
v\to \gu_\a =(\gu_{\a,n} )_{n\in \Z}.
\end{aligned}
$$
for all $n\in\Z$. We first present results on isomorphic inverse
problems  for one-spectra mappings together with norming mappings.

\begin{theorem} \lb{Tisoe+n} Let $\a\in\A$. Then

\no i) The mappings $\m_\a\ts (\pm\gr_\a )$, $\n_\a \ts (\pm\gs_\a
):  \mH\to \gJ^o_\a\ts \ell^2$ are isomorphic. Each of them is a RAB
from $\mH$ onto $\gJ^o_\a\ts \ell^2$ and satisfies
\[
\lb{SpN1}
\begin{aligned}
\m_\a\ts \gr_\a=(\n_\a\ts\gs_\a)\circ\1_- =\big(\m_\a\ts(-\gr_\a)\big)\circ \mU_1 \\
=\big(\n_\a\ts(-\gs_\a)\big)\circ \mU_2 =\big((\cS_\a\m)\ts\gr\big)\ci F_{2\a}^*.
\end{aligned}
\]

\no ii) Each of the mappings $\vr_\a \ts (\pm\gu_\a ) $,
$\t_\a\ts(\pm\gt_\a ): \mH\to\gJ^1_\a\ts \ell^2$ is a RAB from $\mH$
onto $\gJ^1_\a\ts \ell^2$. Moreover, they are isomorphic and satisfy
\[
\lb{SpN2}
\begin{aligned}
\t_\a\ts\gt_\a=(\vr_\a\ts\gu_\a)\circ\1_- =\big(\vr_\a\ts(-\gu_\a)\big)\circ\mU_1
\\
=\big(\t_\a \ts(-\gt_\a)\big)\circ\mU_2 =\big((\cS_\a\t)\times\gt\big)\ci F_{2\a}^*.
\end{aligned}
\]
\end{theorem}

\no {\bf Remark.} This theorem shows that the inverse problems are
isomorphic,  and the corresponding automorphisms are obtained in
explicit form. Each of these inverse problems is a real-analytic
bijection between the space of potentials and the corresponding
spectral data.

We introduce the matrix $J_1=\diag(1,-1)$ and the symmetry subspaces
\[
\lb{def-ha}
\begin{aligned}
\mH_{eo}&=\{v\in\mH:v=J_1\cR v\}, \qqq
\mH_\a&=\{v\in\mH:v=\mU_1v\}.
\end{aligned}
\]
Here even and odd mean even and odd about the midpoint $x=1/2$.
The symmetry space depends on the choice of spectrum.

\begin{corollary}
\label{Ceoa} Let $\a\in\A$. Then each of the mappings $\cS_\a\m,
\cS_\a\n:  \mH_{eo}\to \gJ^o_\a$ is a RAB from $\mH_{eo}$ onto
$\gJ^o_\a$. Moreover, the mappings $\cS_\a\t$, $\cS_\a\vr:
\mH_{eo}\to \gJ^1_\a$ are RABs from $\mH_{eo}$ onto $\gJ^1_\a$. They
satisfy on the space $\mH_{eo}$:
\[
\lb{shifteo}
\cS_\a(\m,\n,\t,\vr)=(\m_\a,\n_\a ,\t_\a,\vr_\a)\circ F_{2\a},
\]
where the mapping $F_{2\a} :\mH_{eo}\to \mH_{\a}$ is a bijection
from $\mH_{eo}$ onto $\mH_{\a}$.
\end{corollary}

\no {\bf Remark.} Thus the mappings $\m_\a$, $\n_\a$, $\t_\a$, and $\vr_\a$ are isomorphic on $\mH_{eo}$.
We are not aware of an analogous result for Schr\"odinger operators on the unit interval $[0,1]$; see
\cite{K25}.

{\bf Replacing mappings.} Recall that the operator $F_{2\a}$ is
defined by \er{DefF_a}. We define  the mapping $\cU_{\a,\s}:
\mH\to\mH$ by
\[
\label{defUa}
\cU_{\a,\s}=F_{2\a}\ci\cU_\s\ci F_{2\a}^*.
\]
 Here $\cU_\s$ denotes the lamplighter mapping introduced in Sec.
2. Its properties will be discussed later. Since $\cU_\s$ is a RAB
from $\mH$ onto itself and $F_{2\a}$ is unitary,  the representation
in \er{defUa} shows that $\cU_{\a,\s}: \mH\to\mH$ is also a RAB from
$\mH$ onto itself.

We now discuss a new type of inverse problem. Let the one-spectra
mapping $v\to \m_\a=(\m_{\a,n})_{n\in \Z}$ be given,    and suppose
we replace some $\m_{\a,n}$ by the eigenvalues $\n_{\a,n}$. Then we
obtain a new mapping $\gc$. For example, we have
$\gc=(...,\m_{\a,1},\n_{\a,2}, \n_{\a,3},\m_{\a,4}, \m_{\a,5},...)$.
{{The question is whether the selected eigenvalues still form an
admissible spectral sequence and, when paired with a complementary
sequence, give a bijective parametrization of the potential. We call
these \emph{spectral replacing mappings}.}} Let
$\Z_1,\Z_2\subset\Z$ be arbitrary subsets. We define replacing
mappings $ v\to\z_\a =(\z_{\a,n})_{n\in\Z}, v\to\phi_\a
=(\phi_{\a,n})_{n\in\Z}$ and their components by
\[
\lb{Repa}
\begin{aligned}
\z_{\a,n}=\left\{\begin{aligned}
&\vr_{\a,n}, \qq \text{if}  &&n\in\Z_1 \\
&\t_{\a,n}, \qq \text{if}  &&n\notin\Z_1
\end{aligned}\right., \qq
\qq
\phi_{\a,n}=\left\{\begin{aligned}
&\n_{\a,n}, \qq \text{if}  &&n\in\Z_2 \\
&\m_{\a,n}, \qq \text{if}  &&n\notin\Z_2
\end{aligned}\right..
\end{aligned}
\]
Define the diagonal sign-change operator $\s=\diag(\s_n)_{n\in\Z}$ by
$$
\s_{2j-1}=\ca  -1 & \text{if} \ j\in \Z_1
\\
1  & \text{if} \ j\notin \Z_1   \ac,\qq
\s_{2j}=\ca  -1 &  \text{if} \ j\in \Z_2
\\
1 & \text{if} \ j\notin \Z_2  \ac.
$$
\begin{corollary} \lb{Creplacing} Let $\a\in\A$.
Define the mappings $\z_\a =(\z_{\a,n})_{n\in\Z}$ and
$\phi_\a =(\phi_{\a,n})_{n\in\Z}$ by \er{Repa}. Then
the mapping $\z_\a \star \phi_\a : \mH\to\gJ_\a $ is a RAB
from $\mH$ onto $\gJ_\a$ and satisfies
\begin{equation}
\label{rma}
(\z_\a \star \phi_\a )\circ \cU_{\a,\s} =\t_\a\star \m_\a.
\end{equation}
\end{corollary}

For two functions $A(v), B(v)$, differentiable on $\mH$, we define the symplectic form
$$
A\wedge B=\int _0^1 (A'(v),B'(v))_0dx,\ \ \ A'(v)={\pa A(v)\/\pa v}.
$$
where the function $ (p,q)_0=p^\top Jq $ for column vectors
$p,q\in\R^2$.  In particular, $A\wedge A=0$. The following canonical
relations identify conjugate spectral coordinates.

\begin{theorem}  \lb{Tcanonicala}
Let $\a\in\A$ and $v\in\mH$. Then for all $n, j\in \Z$,
\[
\lb{cv1a}
\begin{aligned}
\m_{\a,n}(v)\we \m_{\a,j}(v)=0,\qqq \gr_{\a,n}(v)\we \m_{\a,j}(v)=\d_{nj},
\qqq \gr_{\a,n}(v)\we \gr_{\a,j}(v)=0,
\\
\n_{\a,n}(v)\we \n_{\a,j}(v)=0,\qqq \gs_{\a,n}(v)\we \n_{\a,j}(v)=\d_{nj}, \qqq
\gs_{\a,n}(v)\we \gs_{\a,j}(v)=0,
\\
\t_{\a,n}(v)\we \t_{\a,j}(v)=0,\qqq \gt_{\a,n}(v)\we \t_{\a,j}(v)=\d_{nj}, \qqq
\gt_{\a,n}(v)\we \gt_{\a,j}(v)=0,
\\
\vr_{\a,n}(v)\we \vr_{\a,j}(v)=0,\qqq \gu_{\a,n}(v)\we \vr_{\a,j}(v)=\d_{nj}, \qqq
\gu_{\a,n}(v)\we \gu_{\a,j}(v)=0.
\end{aligned}
\]
{{Moreover, each of the following four families  $\{\m'_{\a,n},
\gr'_{\a,n}\}_{n\in \Z}$, \ $\{\n'_{\a,n},$ $ \gs'_{\a,n}\}_{n\in
\Z}$,  $\{\t'_{\a,n}, \gt'_{\a,n}\}_{n\in \Z}$, \ $\{\vr'_{\a,n},
\gu'_{\a,n}\}_{n\in \Z}$ is a Riesz basis of $\mH$}},
 where $\m'_{\a,n}={\pa\m_{\a,n}\/\pa
v},\gr'_{\a,n}={\pa\gr_{\a,n} \/\pa v}$,... for brevity.
\end{theorem}

\no  Note that such canonical relations were established for
Schr\"odinger operators in \cite{FM76,PT87}, while those for ZS
operators under basic boundary conditions were obtained in
\cite{K01,K05,KZ26}.

\subsection{Short review}
We briefly review known results in inverse spectral theory for
differential  operators on the unit interval. Inverse
Sturm--Liouville  problems have been studied by many authors,
including Borg \cite{B46}, Gel'fand-Levitan \cite{GL51}, Krein
\cite{Kr51}, and Marchenko \cite{M50}. But in all of them there was
a gap between necessary and sufficient conditions.
Marchenko and Ostrovski \cite{MO75} closed this gap and solved inverse problems for
the two-spectra  mapping including a characterization, see also
Levitan and Gasymov \cite{LG64} about it. These results form part of
the background to future investigations. Trubowitz and his coauthors
\cite{IT83, IMT84, PT87, T77} proposed analytic approach based on
 properties of the map from potentials to spectral data, as
well as explicit transformations that alter only {\it finitely} many
spectral parameters (eigenvalues and norming constants). This
approach was developed further in  and applied to other inverse
problems with purely discrete spectra: (a)  various Sturm--Liouville
problems \cite{CM93, CM93a, GT87, KC09}, (b)
 perturbed harmonic oscillators, see \cite{MT81, CKK04, CK07};
(c) Sturm--Liouville problems with matrix-valued potentials subject
to Dirichlet boundary conditions on the unit interval, see
\cite{CK06q, CK09}; (d) Jacobi matrices, see \cite{KL27} and the
references therein.

Marchenko and Ostrovski \cite{MO75} also solved inverse periodic
problems in terms of the conformal mappings. Dubrovin \cite{Du75}
obtained the famous equation to recover the periodic potential via
spectral data, see also \cite{T77}. An analytic approach from
\cite{PT87} was applied by Kargaev and Korotyaev \cite{KK97} for even
periodic potentials, and for general case by Korotyaev \cite{K97, K99}, including periodic distributions \cite{K03}. In this
case a priori estimates of potentials in terms of gap lengths (or
other spectral parameters) are important in the proof. Two-sided
estimates for the potential in terms of gap lengths  (or parameters
of the Lyapunov function) were obtained in \cite{KK00,
K06,K05,K01,K00a,K96,K98} using the conformal mapping theory
associated with quasimomentum. Various properties of the relevant
conformal mappings were established in \cite{CK06, KK97, KK95, K97m,
MO75}. We recall that Hilbert \cite{H09} obtained the first result
on conformal mappings from multiply connected domains onto domains
with parallel slits.

 There are lot of results about periodic systems, Lyapunov,
Gel'fand-Lidskii \cite{GL55}, Gesztesy and coauthors \cite{CG02, C00}, Krein \cite{Kr55},  Chelkak-Korotyaev \cite{CK06q},
Korotyaev \cite{K08, K10} and see  references therein.
 The matrix potentials pose new problems. The Riemann
surfaces are multi-sheeted, have complicated structures and their
branch points can be real and complex. The spectral bands can have
various multiplicity and their end points are branch points or
periodic eigenvalues. Their control is not simple. Until now the
inverse problem on the finite interval  is solved only for the
Schr\"odinger operator with matrix-valued potential  \cite{CK06q, CK09} (including characterization) and for first order systems
\cite{M99, M15} (only uniqueness). First order periodic
systems with $3\ts 3$ matrix valued potentials are discussed in
\cite{K25q}.
These results are important to study cubic periodic vector NLS equations.

Now we discuss  inverse problems for ZS operators, which have also
been studied extensively. For systems on the unit interval subject
to boundary conditions, we refer  to
\cite{AHM05,A14,A93,DK00,GJ75,H01,LS91,MT02,S70}; for systems on the
circle, see \cite{BG93,K05,K01,K96,KM22,KM21,LS91,GG93} and the
references therein. Uniqueness and reconstruction results for these
inverse problems are given in \cite{CG02, GL66, LS91}; see also
\cite{Ge02,H01,RG01, MT02}. Misura obtained the complete results for
two spectra mapping \cite{Mi78,Mi79}. Misura also solved inverse
periodic problems in terms of in terms of the conformal mappings,
\cite{Mi78,Mi79}. In fact she extends  results and the proof of
Machenko-Ostrovki for Hill operators \cite{MO75} fo the Dirac
operators. After long period Albeverio and his coauthors
\cite{AHM05} solved the characterization problems for the
two-spectra mappings $\vr\star\n$ and the mapping $\vr$ plus the
corresponding normalazing constants. Using analytic approach
Korotyaev \cite{K01} reproved the results of Misura for periodic
case. Moreover, he solved the inverse periodic problem  in terms gap lengths in \cite{K05}.
The estimates of potentials in terms of gap lengths (and other parameters) were obtained in \cite{K01,K05, K96}.

Finally, we discuss isomorphic inverse problems. Korotyaev
\cite{K25} described isomorphic inverse problems for
Sturm--Liouville operators on the unit interval and on the circle.
Related problems for operators in impedance form are the subject of
the work in preparation \cite{K27}. Isomorphic inverse problems are
discussed for Jacobi operators in \cite{KL27}, and for ZS operators
on the unit interval and the circle in \cite{KZ26}. These works
introduce and use four-spectra mappings and lamplighter mappings,
which serve as important tools for the study of isomorphic inverse
problems.

\section {Preliminary results \lb{Sec2}}
\setcounter{equation}{0}

\subsection{Basic problems for $\a=0$}
 For $\a=0$, problems \er{sgc1}-\er{sgc4}
reduce to four corresponding basic problems. In this section, we
introduce the notation and the results for the basic problems that
will be used in the gauge reduction. These results are recalled from
\cite{KZ26}; they are inputs to the arguments below, rather than new
assertions proved here.

The basic boundary value problems are the ZS systems
$$
Jy'+Vy=\l y,
$$
subject to the boundary conditions \er{1.1'}-\er{1.4'}. These four
conditions are the special cases of \er{sgc1}-\er{sgc4}  when
$\a=0$, namely the Dirichlet and Neumann boundary conditions:
\[
\lb{1.1'}
y_1(0)= y_1(1)=0,
\]
\[
\lb{1.2'}
y_2(0)=y_2(1)=0,
\]
and the so-called mixed boundary conditions:
\[
\lb{1.3'}
y_1(0)=y_2(1)=0,
\]
\[
\lb{1.4'}
y_2(0)=y_1(1)=0,
\]
respectively. Let $\mu_n$, $\nu_n$, $\tau_n$ and $\vr_n$ ($n\in \Z$)
denote the eigenvalues corresponding to these four problems. All
these eigenvalues are simple and satisfy
$$
\begin{aligned}
& ......<\ol {\tau_1, \vr_1}<\ol {\mu_1,\nu_1}<\ol {\tau_2, \vr_2}< \ol
{\mu_2,\nu_2}<...,
\\
& \nu_n,\mu_n=\mu_n^o+o(1),\qq \tau_n, \vr_n=\tau_n^o+o(1)\qq \as \ n\to
\pm \infty,
\end{aligned}
$$
where the ordered-pair convention is the same  as in Sect. 1, and
$\nu_n^o=\mu_n^o=\pi n$ and $\tau_n^o=\vr_n^o=\pi (n-\frac{1}{2}),
n\in \Z$ are the corresponding unperturbed eigenvalues.

When $\a=0$, the families of real strictly increasing sequences
$\gJ^o_\a$, $\gJ^1_\a$, $\gJ_\a$ reduce to $\gJ^o$, $\gJ^1$, $\gJ$,
respectively,
$$
\begin{aligned}
&\gJ^o=\Big\{ s\!=\!(s_n)_{n\in \Z}: ....<s_1<s_2<....., \qq
s_n\!=\!\pi n\!+\breve{s}_{n\,}   ,\qqq \breve{s}=(
\breve{s}_{n\,})_{n\in \Z}\!\in\!\ell^2     \Big\},
\\
& \tes   \gJ^1=\Big\{ s\!=\!(s_n)_{n\in \Z}: ....<s_1<s_2<....., \qq
s_n\!=\!\!\pi (n-\frac{1}{2})\!+\breve{s}_{n\,},\qqq \breve{s}=(
\breve{s}_{n\,})_{n\in \Z}\!\in\!\ell^2     \Big\},
\\
&\gJ=\Big\{ t\!=\!(t_n)_{n\in \Z}: ...<t_1<t_2<....., \qq
t_n\!=\frac{\pi n}{2}+\breve{t}_{n\,},\qqq \breve{t}=(
\breve{t}_{n\,})_{n\in \Z}\!\in\!\ell^2 \Big\}.
\end{aligned}
$$

\subsection{Results for $\a=0$}
We introduce the four-spectra mapping $\gf: \mH\to \ell^2$ from
\cite{KZ26}. It is defined via the basic eigenvalues and takes the
following form:
\begin{equation}
\label{Def4sm}
\tes
v\to \gf(v)=(\gf_n(v))_{n\in \Z},\qq \gf_{2n-1}=\frac{1}{2}(\vr_n
-\tau_n),\qqq  \gf_{2n}=\frac{1}{2}(\nu_n -\mu_n), \ \ \ \ \ n\in \Z.
\end{equation}

We recall a theorem proved in \cite{KZ26}.

\begin{theorem} \lb{L4sm}
The four-spectra mapping $\gf: \mH\to \ell^2$  is a RAB from $\mH$
onto $\ell^2$ and satisfies
\begin{equation}
\label{egf}
\begin{aligned}
 \tes
 \frac{1}{\sqrt 2}\|\gf(v)\|\leq \|v\|\leq 2\|\gf(v)\|\big(1+\|\gf(v)\|\big),\qqq \forall\ v\in\mH.
\end{aligned}
\end{equation}
\end{theorem}

\no Let $\mD$ be the set of all $v\in\mH$ such that the sequence
$(\gf_n(v))_{n\in\Z}$ is finitely supported, i.e., $\gf_n(v)=0$ for
all sufficiently large $|n|$. Since the four-spectra mapping
$\gf:\mH\to\ell^2$ is a RAB from $\mH$ onto $\ell^2(\Z)$, the set
$\mD$ is dense in $\mH$. Let $\gS$ be the set of all diagonal
operators  $\sigma=\diag (\ldots,\s_1,\s_2,\ldots)$ on $\ell^2$,
written briefly as $\s=(\s_j)_{j\in \Z}$, where $\s_j\in\{\pm1\}$
for every $j\in\Z$. The set $\gS$ defines the so-called lamplighter
group; see \cite{CM17}. For each $\s\in\gS$, with $\gf$ given by
\er{Def4sm}, we define the lamplighter mapping $\cU_\s$ by
$$
 \cU_\s=\gf^{-1} \ci(\s \gf): \mH\to \mH.
$$
Note that $\cU_\s:\mD\to\mD$ for every $\s\in\gS$. We define a
reflection  $\cR$  and reflection-type operators $\cR_1, \cR_2$
acting on $\mH$ by
$$
(\cR v)(x)= v(1-x),\qq x\in[0,1], \qq {\rm and} \qqq
\cR_1=J_1\cR,  \qqq  \cR_2=-\cR_1.
$$

\begin{theorem} \label{Llamp}
i) Each mapping $\cU_\s:\mH\to\mH, \s\in\gS$ is a RAB of $\mH$ onto
itself and satisfies
$$
\begin{aligned}
\cU_\s=\cU_\s^{-1}, \qqq
\cU_\s\ci \cU_{\s'}=\cU_{\s\s'}=\cU_{\s'}\ci
\cU_{\s} \qqq \forall \ \s, \s'\in \gS,
\end{aligned}
$$
$$
\|\cU_\s(v)\|=\|v\|\qqq \forall \ v\in\mH.
$$
\no ii) Let $(\l_n^\pm)_{n\in\Z}$ denote  the periodic and
antiperiodic  eigenvalues of the $1$-periodic extension of $v$,
equivalently the 2-periodic eigenvalues. They are invariant under
each $\cU_\s, \s\in \gS$, i.e.,
$$
(\l_{n}^\pm)_{n\in\Z}=(\l_{n}^\pm)_{n\in\Z}\ci \cU_\s
\qq \forall \ \s\in \gS,
$$
and the mappings $\cU_\s$ for specific $\s\in\gS$ have the forms:
$$
\begin{aligned}
\1_-=\cU_\s,\qq \where \ \s=\1_-,\qqq
\ca \cR_1=\cU_\s,\qq \where \qq \s_{n}= (-1)^n \qqq
\forall \ n\in \Z \\
\cR_2=\cU_\s,\qq \where \qq \s_{n}= -(-1)^n \qqq
\forall \ n\in \Z\ac.
\end{aligned}
$$

\no iii) If  $\s=(\s_j)_{j\in\Z}\in\gS$, then for each $n,j\in \Z$ we have
$$
{\rm if} \ n=2j-1 \Rightarrow \ca (\tau_j, \vr_j)=(\varrho_j,\t_j)\ci
\cU_\s, \ \ {\rm if}\qq  \s_n=-1 \\
(\tau_j, \vr_j)=(\t_j,\varrho_j)\ci \cU_\s, \ \ {\rm if} \qq
\s_n=1 \ac,
$$
$$
{\rm if} \ n=2j \Rightarrow \ca (\m_j, \n_j)=(\n_j,\m_j)\ci \cU_\s, \ \ {\rm if}\qq  \s_n=-1 \\
(\m_j, \n_j)=(\m_j,\n_j)\ci \cU_\s \ \ {\rm if} \qq \s_n=1 \ac .
$$
\end{theorem}

We define the one-spectra mappings $\mu, \nu: \mH\to\gJ^o$ and
$\t,\vr: \mH\to\gJ^1$ by
$$
\begin{aligned}
v\to \mu=(\mu_n)_{n\in \Z},\ \qq  v\to \nu=(\nu_n)_{n\in \Z},
\qq
 v\to \t=(\tau_n)_{n\in \Z},\ \qq v\to \vr=(\vr_n)_{n\in \Z}.
 \end{aligned}
$$
Combining two one-spectra mappings $\t=(\tau_n)_{n\in \Z}\in \gJ^1$
and  $\mu=(\mu_n)_{n\in \Z} \in \gJ^o$, whose sequences are strictly
increasing and interlaced, we construct a two-spectra mapping
$\t\star\m$ from $\mH$ into $\gJ$ by
$$
v\to \tau \star \mu=(....,\tau_1,\mu_1,\tau_2,\mu_2,....).
$$
We introduce the real norming constants $\gr_{n}, \gs_{n}$,
$\gt_{n},  \gu_{n}$, $n\in\Z$ and the corresponding {\it norming }
mappings by
$$
\begin{aligned}
& \gr_n=-\ln |\vp_2(1,\m_n)|, \qq \gs_n=-\ln |\vt_1(1,\n_n)|, \qq
  \gt_n=-\ln |\vp_1(1,\t_n)|, \qq  \gu_n=-\ln |\vt_2(1,\vr_n)|,
 \\
& v\to \gr=(\gr_n )_{n\in \Z},\qq  v\to \gs=(\gs_n )_{n\in \Z},\qq
  v\to \gt=(\gt_n )_{n\in \Z},\qq  v\to \gu=(\gu_n )_{n\in \Z}.
\end{aligned}
$$
We recall the isomorphism for the basic boundary value problems established in \cite{KZ26}.

\begin{theorem} \lb{Liso}
\no i) All two-spectra mappings $\t\star \m, \vr\star \n$, $\vr\star
\m$ and $\t\star \n$ from $\mH$ to $\gJ$  are isomorphic, each of
them is a RAB from $\mH$ onto $\gJ$ and they satisfy
$$
\begin{aligned}
\t\star \m=(\vr\star \n)\circ\1_-=(\vr\star \m)\circ \cR_1=(\t\star
\n)\circ \cR_2.
\end{aligned}
$$
\no ii)  Each of the mappings $\m\ts (\pm\gr)$, $\n \ts (\pm\gs)$
from $\mH$ to $\gJ^o\ts \ell^2$ is a RAB from $\mH$ onto
$\gJ^o\times\ell^2$. Moreover, they are isomorphic and satisfy
$$
\begin{aligned}
\m\ts \gr =(\n \ts \gs )\circ\1_-
=(\m \ts (-\gr) )\circ \cR_1=(\n \ts (-\gs) )\circ \cR_2.
\end{aligned}
$$
\no iii)  Each of the mappings $\t\ts (\pm\gt), \vr \ts (\pm\gu ) $,
from $\mH$ to $\gJ^1\ts \ell^2$ is a RAB from $\mH$ onto $\gJ^1\times\ell^2$.
Moreover, they are isomorphic and satisfy
$$
\begin{aligned}
\t\ts \gt =(\vr \ts \gu )\circ\1_-
=(\vr \ts (-\gu) )\circ \cR_1=(\t \ts (-\gt) )\circ \cR_2.
\end{aligned}
$$
\end{theorem}

The following one-spectra result concerns even--odd potentials
with respect to $x=1/2$.

\begin{theorem}
\label{L1seo} Each of the mappings $\cS\m$, $\cS\n$, $\t$, and $\vr$
from $\mH_{eo}$ to $\gJ^1$ is a RAB from $\mH_{eo}$ onto $\gJ^1$. On
$\mH_{eo}$, they satisfy
$$
\cS\m=(\cS\n)\ci \cF_o=\t\circ \gF=\vr\ci \cF_o\circ \gF,
$$
where the operator $\gF:\mH_{eo}\to\mH_{oe}$ is a bijection from
$\mH_{eo}$ onto $\mH_{oe}$.
 \end{theorem}

Following Marchenko \cite{M50}, we define the real normalizing
constants $\ga_n$, $\gb_n$, $\gc_n$, and $\gd_n$, $n\in\Z$, together
with the associated {\it normalizing} mappings as follows
\begin{equation}
\label{mnc}
\begin{aligned}
&e^{-\ga_n}=\|\vp(\cdot,\mu_n)\|^2,\qq
e^{-\gb_n}=\|\vt(\cdot,\nu_n)\|^2, \qq
e^{-\gc_n}=\|\vp(\cdot,\t_n)\|^2,\qq
e^{-\gd_n}= \|\vt(\cdot,\vr_n)\|^2,
\\
& v\to\ga=(\ga_n )_{n\in \Z}, \qq   v\to\gb =(\gb_n)_{n\in \Z},\qq
v\to\gc=(\gc_n )_{n\in \Z}, \qq   v\to\gd=(\gd_n)_{n\in \Z}.
\end{aligned}
\end{equation}
Note that $\ga(0)=\gb(0)=\gc(0)=\gd(0)=0$. We discuss inverse
problems for one-spectra mappings and normalizing  mappings,
following \cite{KZ26}.

\begin{theorem} \lb{Le+nl}
i) Each of the mappings $\m\ts \ga$ and $\n \ts \gb$, defined by the
basic spectral mappings above and \er{mnc}, is a bijection from
$\mH$ onto $\gJ^o\ts \ell^2$ and satisfies
$$
\m\ts \ga=(\n \ts \gb)\circ \1_-.
$$
\no ii) Each of the mappings $\t\ts \gc$ and $\vr \ts \gd$, defined
by the basic spectral mappings above and \er{mnc}, is a bijection
from $\mH$ onto $\gJ^1\ts \ell^2$ and satisfies
$$
\t\ts \gc=(\vr \ts \gd)\circ \1_-.
$$
\end{theorem}

The following canonical relations for the basic eigenvalues were proved in \cite{KZ26}.

\begin{theorem}  \lb{Lcanonical}
Let $v\in\mH$. Then for any $n, j\in \Z$, the following identities
hold:
$$
\begin{aligned}
\m_n(v)\we \m_j(v)=0,\qqq \gr_n(v)\we \m_j(v)=\d_{nj},
\qqq \gr_n(v)\we \gr_j(v)=0, \\
\n_n(v)\we \n_j(v)=0,\qqq \gs_{n}(v)\we \n_j(v)=\d_{nj}, \qqq
\gs_{n}(v)\we \gs_{j}(v)=0, \\
\t_n(v)\we \t_j(v)=0,\qqq \gt_{n}(v)\we \t_j(v)=\d_{nj}, \qqq
\gt_{n}(v)\we \gt_{j}(v)=0, \\
\vr_n(v)\we \vr_j(v)=0,\qqq \gu_{n}(v)\we \vr_j(v)=\d_{nj}, \qqq
\gu_{n}(v)\we \gu_{j}(v)=0.
\end{aligned}
$$
Moreover, each of the four gradient families
$$
\{\m_n',\gr_n'\}_{n\in\Z},\quad
\{\n_n',\gs_n'\}_{n\in\Z},\quad
\{\t_n',\gt_n'\}_{n\in\Z},\quad
\{\vr_n',\gu_n'\}_{n\in\Z}
$$
is a Riesz basis of $\mH$,  where $\m'_{\a,n}={\pa\m_{\a,n}\/\pa
v},\gr'_{\a,n}={\pa\gr_{\a,n} \/\pa v}$,... for brevity.

\end{theorem}

Let $\Y_1$ and
$\Y_2$ be subsets of $\Z$. We define the replacing mappings
$v\to\z=(\z_n)_{n\in\Z}$ and $v\to\phi=(\phi_n)_{n\in\Z}$ by
\[
\lb{Rep}
\begin{aligned}
\z_n=\left\{\begin{aligned}
&\vr_n, \qq \text{if}  &&n\in\Y_1 \\
&\t_n, \qq \text{if}  &&n\notin\Y_1
\end{aligned}\right., \qq
\qq
\phi_n=\left\{\begin{aligned}
&\n_n, \qq \text{if}  &&n\in\Y_2 \\
&\m_n, \qq \text{if}  &&n\notin\Y_2
\end{aligned}\right..
\end{aligned}
\]
Define the diagonal sign-change operator $\s=\diag(\s_n)_{n\in\Z}$ by
$$
\s_{2j-1}=\ca  -1  &  \text{if} \ j\in \Y_1
\\
1  & \text{if} \ j\notin\Y_1 \ac, \qqq
\s_{2j}=\ca  -1   &  \text{if} \ j\in \Y_2
\\
1 & \text{if} \ j\notin \Y_2 \ac.
$$

\begin{theorem} \lb{Lreplacing}
Let the mappings $\z=(\z_n)_{n\in\Z}$ and $\phi=(\phi_n)_{n\in\Z}$
be defined by \er{Rep}. Then $\z\star\phi:\mH\to\gJ$ is a RAB from
$\mH$ onto $\gJ$ and satisfies
$$
(\z\star \phi)\circ \cU_\s=\t\star \m.
$$
\end{theorem}

\subsection{Generic boundary conditions}
We consider the following generic boundary value problem for the
vector-valued function $f=(f_1, f_2)^\top$
\[
\lb{gbp} Jf'+Vf=\l f, \qqq  \left\{\begin{aligned}
f_1(0)c_\b-f_2(0)s_\b=0  \\ f_1(1)c_\a -f_2(1)s_\a
=0\end{aligned}\right. .
\]
We show that this generic problem is unitarily equivalent to the
problem  \er{sgc1}, with a transformed potential and a shifted
spectral parameter. Similarly, there exist corresponding generic
problems that are unitarily equivalent to the problems
\er{sgc2}--\er{sgc4}, respectively.

For $\b\in\A$, let $\gF_{\b}, \gF_{2\b}: \mH\to\mH$ denote the
unitary  multiplication operators defined by
\[
\lb{defFb} \gF_\b f=e^{\b(x-1) J} f, \qq \text{and} \qq \gF_{2\b}
f=e^{2\b(x-1) J} f,  \qqq \gF_{2\b}=\gF_\b^2.
\]
A direct calculation yields
$$
\gF_\b \big(J\tes\frac{d}{dx}+V\big)\gF_\b ^*=J\tes\frac{d}{dx}+\b\1_2+V_\b,
$$
where
\[
\lb{Vb}
V_\b=\gF_\b V \gF_\b^*=\gF_{2\b} V=\ma u_1 & u_2\\ u_2 &-u_1\am,
\]
and
$$
u=\gF_{2\b} v=\ma v_1c+v_2s \\ -v_1s+v_2c \am, \qqq
\ca c=\cos\big(2\b(x-1)\big) \\ s=\sin\big(2\b(x-1)\big) \ac .
$$
The following lemma describes how the general boundary value problem
transforms under the unitary transformation $\gF_{\b}$.

\begin{lemma} \lb{GTsG}
Let $\a, \b\in\A$ and $v\in\mH$. Let $\gF_\b(x)$ be defined by
\er{defFb}.

\no i) A function \(f\) solves the boundary value problem
$$
\tes Jf'+Vf=\l f, \qq x\in(0,1),\qqq \left\{\begin{aligned} f_1(0)c_\b-f_2(0)s_\b=0
\\ f_1(1)c_\a -f_2(1)s_\a =0\end{aligned}\right. ,
$$
if and only if $y=\gF_\b f$ solves
$$
Jy'+V_\b y=(\l-\b)y,\qq x\in(0,1),\qqq
\left\{\begin{aligned}&y_1(0)=0\\ &y_1(1)c_\a -y_2(1)s_\a =0\end{aligned}\right. ,
$$
where $V_\b$ is given by \er{Vb}.

\no ii) The fundamental matrix solution $Y(x,\l,u)$, where $u=\gF_{2\b} v$, satisfies
\[
\lb{funds} Y(1,\l-\b,u)= \ma c_\b\vt_1-s_\b\vp_1 &
s_\b\vt_1+c_\b\vp_1 \\ c_\b\vt_2-s_\b\vp_2  & s_\b\vt_2+c_\b\vp_2
\am (1,\l,v), \qqq Y(0,\cdot,u)=\1_2.
\]
\end{lemma}

\no {\bf Remark.} This lemma shows that the generic boundary
condition  $f_1(0)c_\b-f_2(0)s_\b=0$ simplifies to the Dirichlet
condition $y_1(0)=0$, while the boundary condition at $x=1$ retains
its original structure.

\no {\bf{Proof.}} i) Let $\gF_\b(x)=e^{\b(x-1)J}$.
The gauge transformation of system \er{gbp} yields
\[
\lb{sgbp}
\begin{aligned}
\gF_\b(J\frac{d}{dx}+V)\gF_\b^*y
=Jy'+V_\by=(\l-\b)y,\qqq
V_\b=\gF_\b V\gF_\b^*.
\end{aligned}
\]
By the standard existence and uniqueness theorem for linear ODEs, we obtain
$$
y(x,\l-\b,u)=\gF_\b(x) f(x,\l,v), \qqq \text{where} \qqq u=\gF_{2\b}v.
$$
A direct computation using the preceding identity gives
$$
\begin{aligned}
&y(0,\l-\b,u)=\gF_\b(0)f(0,\l,v)=\ma c_\b f_1-s_\b f_2\\ s_\b f_1 +c_\b f_2 \am(0,\l,v), \\
&y(1,\l-\b,u)=\gF_\b(1)f(1,\l,v)=\ma f_1  \\ f_2  \am(1,\l,v).
\end{aligned}
$$
We therefore conclude that the transformed solution satisfies
$$
y_1(0)=0 \qqq \text{and} \qqq y_1(1)c_\a -y_2(1)s_\a =0.
$$
This establishes part i).

ii) By \er{sgbp}, we deduce that the fundamental matrix solution $Y$ satisfies
$$
Y(x,\l-\b,u)=\gF_\b(x) Y(x,\l,v)e^{\b J}, \qqq Y(0,\cdot,u)=\1_2,
$$
which implies \er{funds}.           \BBox

\section{Gauge transformations and their properties \lb{Sec3}}
\setcounter{equation}{0}

In this section, we study the properties of ZS systems under
different  unitary transformations. We describe the behavior of
spectral data, eigenfunctions, and other related quantities under
such transformations.

Using the unitary operator $F_\a(x)=e^{\a xJ}$ acting on $\mH$ ,  we
obtain that
$$
F_\a\big(J\frac{d}{dx}+V\big)F_\a^*=J\frac{d}{dx}+\a\1_2+V_\a,
$$
where
$$
V_\a=F_\a VF_\a^*=F_{2\a} V=\ma u_1 & u_2 \\ u_2 & -u_1\am,
$$
and
$$
u=\ma u_1 \\ u_2\am=F_{2\a} v=\ma c_{2\a x}v_1+s_{2\a x}v_2
\\
-s_{2\a x}v_1+c_{2\a x}v_2 \am.
$$
We investigate the correspondence between the spectra of the
original operators and those of the transformed ones.

\begin{lemma}
\lb{LrotationF_a} Let $\a\in\A$ and $v\in\mH$. Let $F_{2\a}$ and
$\cS_\a$ be defined by \er{DefF_a} and \er{DefShiftO}, respectively.
Then

\no i) The shifted operator $\cS_\a:\gJ\to\gJ_\a$ is a real-analytic
isomorphism from $\gJ$ onto $\gJ_\a$. The same holds for the
corresponding translations of $\gJ^o$ and $\gJ^1$.

\no ii) The fundamental matrix solution $Y$ has the following form
$$
Y(x,\l-\a, u) =\ma c_{\a x} \vt_1+s_{\a x} \vt_2 & c_{\a x}
\vp_1+s_{\a x} \vp_2  \\ -s_{\a x} \vt_1+c_{\a x} \vt_2 & -s_{\a x}
\vp_1+c_{\a x} \vp_2\am(x,\l,v), \qq u=F_{2\a} v.
$$

\no iii) The one-spectra mappings $\m_\a,\n_\a ,\t_\a, \vr_\a$, and
the norming mappings $\gr_\a$, $\gs_\a$, $\gt_\a$, $\gu_\a$, defined
by \er{DefNorming}, satisfy
\[
\lb{cA}
(\m_\a,\n_\a ,\t_\a,\vr_\a)\circ F_{2\a} =\cS_\a(\m,\n,\t,\vr).
\]
\[
\lb{cA1}
(\gr_\a,\gs_\a,\gt_\a,\gu_\a)\ci F_{2\a}=(\gr,\gs,\gt,\gu).
\]
\end{lemma}

\no {\bf{Proof.}} i) We first prove that $\cS_\a$ is injective.
Suppose that $\cS_\a a=\cS_\a b$ for two sequences
$a=(a_n)_{n\in\Z}$ and $b=(b_n)_{n\in\Z}$ in $\gJ$. By the
definition of $\cS_\a$, we have $a_n=b_n$ for all $n\in\Z$, and
hence $a=b$. To prove a surjection, let $b=(b_n)_{n\in\Z}\in\gJ_\a$
and define $a=(a_n)_{n\in\Z}$ by $a_n=\tes b_n+\a$. Then $a\in\gJ$,
and $\cS_\a a=b$ by construction. Thus $\cS_\a$ is bijective. In
standard coordinates it and its inverse are the identity, so both
are real analytic. The same argument applies to $\gJ^o$ and $\gJ^1$.

ii) Let $F_\a=e^{\a x J}$. The gauge transformation for the system \er{eqZS} gives
\[
\lb{equaFa}
\begin{aligned}
F_\a\big(J\tes\frac{d}{dx}+V\big)F_\a^*f
=Jf'+V_\a f&=(\l-\a)f,
\end{aligned}
\]
where $V_\a=F_\a VF_\a^*=F_{2\a} V$. Let $u=F_{2\a} v$ and
$\z=\l-\a$. By standard ODE theory, we deduce that the fundamental
matrix solution $Y$ of equation \er{equaFa} satisfies
\[
\lb{3.1}
Y(x,\z,u)=F_\a(x) Y(x,\l,v)=\ma c_{\a x} \vt_1+s_{\a x}
\vt_2 & c_{\a x} \vp_1+s_{\a x} \vp_2 \\
-s_{\a x} \vt_1+c_{\a x} \vt_2 & -s_{\a x} \vp_1+c_{\a x} \vp_2\am(x,\l,v).
\]
Recall that, for any $n\in\Z$, the quantities $\m_{\a,n}$,
$\n_{\a,n},\t_{\a,n},\vr_{\a,n}$ denote the eigenvalues of problems
\er{sgc1}-\er{sgc4}, respectively. Then using identity \er{3.1}, we
obtain
$$
\begin{aligned}
\big(c_\a  \vp_1-s_\a \vp_2\big)(1,\z,u)=\vp_1(1,\l,v), \qqq
&\big(s_\a  \vp_1+c_\a \vp_2\big)(1,\z,u)=\vp_2(1,\l,v),
\\
\big(c_\a  \vt_1-s_\a \vt_2\big)(1,\z,u)=\vt_1(1,\l,v), \qqq
&\big(s_\a  \vt_1+c_\a \vt_2\big)(1,\z,u)=\vt_2(1,\l,v),
\end{aligned}
$$
which yields
$$
\begin{aligned}
&\m_{\a,n}(u)=\m_n(v)-\a, \qqq \n_{\a,n}(u)=\n_n(v)-\a, \\
&\t_{\a,n}(u)=\t_n(v)-\a, \qqq \vr_{\a,n}(u)=\vr_n(v)-\a,
\end{aligned}
$$
for all $n\in\Z$. This proves \er{cA}. By the definition of  norming
constants, we obtain for all $n\in\Z$:
$$
\begin{aligned}
&\gr_{\a,n}(u)=-\ln\Big|\frac{\vp_2(1,\m_{\a,n},u)}{c_\a}\Big|
=-\ln|\vp_{2}(1,\m_n,v)|
=\gr_n(v), \\
&\gs_{\a,n}(u)=-\ln\Big|\frac{\vt_1(1,\n_{\a,n},u)}{c_\a}\Big|
=-\ln|\vt_{1}(1,\n_n,v)|
=\gs_n(v), \\
&\gt_{\a,n}(u)=-\ln\Big|\frac{\vp_1(1,\t_{\a,n},u)}{c_\a}\Big|
=-\ln|\vp_{1}(1,\t_n,v)|
=\gt_n(v), \\
&\gu_{\a,n}(u)=-\ln\Big|\frac{\vt_2(1,\vr_{\a,n},u)}{c_\a}\Big|
=-\ln|\vt_{2}
(1,\vr_n,v)|=\gu_n(v).
\end{aligned}$$
These identities imply the identity \er{cA1}.       \BBox

\begin{lemma} \lb{LrotationF_12}
Let $\a\in\A$ and $\mU_1, \mU_2$ be defined by \er{DefF_a},
\er{DefF_12}. Then the one-spectra mappings $\m_\a$, $\n_\a $,
$\t_\a$, $\vr_\a$ and the norming mappings $\gr_\a $, $\gs_\a $,
$\gt_\a $, $\gu_\a $ satisfy
\[
\lb{RelaE}
\begin{aligned}
(\m_\a,\n_\a,\t_\a,\vr_\a)
=(\n_\a,\m_\a,\vr_\a,\t_\a)\circ\1_-
=(\m_\a,\n_\a,\vr_\a,\t_\a)\circ\mU_1
=(\n_\a,\m_\a,\t_\a,\vr_\a)\circ\mU_2,
\end{aligned}
\]
and
\[
\lb{RelaN}
\begin{aligned}
(\gr_\a,\gs_\a,\gt_\a,\gu_\a)=
(\gs_\a, \gr_\a, \gu_\a, \gt_\a)\circ\1_-
=-(\gr_\a,\gs_\a,\gu_\a,\gt_\a)\circ\mU_1
=-(\gs_\a,\gr_\a,\gt_\a,\gu_\a)\circ\mU_2.
\end{aligned}
\]
\end{lemma}

\no {\bf Remark.} Note that $-(\gr_\a,\gs_\a,\gu_\a,\gt_\a)
=(-\gr_\a,-\gs_\a,-\gu_\a,-\gt_\a)$.

\no {\bf Proof.} We prove the identities \er{RelaE}. The gauge
transformation for the system \er{eqZS} gives
$$
\lb{new}
\begin{aligned}
\tes J\big(J\frac{d}{dx}+V\big)J^*f=\tes J f'-Vf=\l f,
\end{aligned}
$$
since $JV=-VJ$. From standard results in ODE theory, we deduce that
the fundamental matrix solution $Y$ of the preceding equation
satisfies
$$
Y(1,\l,p)=JY(1,\l,v)J^*=\ma \vp_2& -\vt_2\\ -\vp_1&\vt_1\am(1,\l,v), \qqq p=\1_-v.
$$
Expressing this matrix identity in terms of the entries of \(Y\), we get
\[
\lb{RelaF}
\begin{aligned}
(c_\a\vp_1-s_\a\vp_2)(1,\l,p)=-(c_\a\vt_2+s_\a\vt_1)(1,\l,v),
\\
(s_\a\vt_1+c_\a\vt_2)(1,\l,p)=-(c_\a\vp_1-s_\a\vp_2)(1,\l,v),
\\
(s_\a\vp_1+c_\a\vp_2)(1,\l,p)=(c_\a\vt_1-s_\a\vt_2)(1,\l,v),
\\
(c_\a\vt_1-s_\a\vt_2)(1,\l,p)=(s_\a\vp_1+c_\a\vp_2)(1,\l,v).
\end{aligned}
\]
Consequently, the associated eigenvalues satisfy, for every \(n\in\Z\),
\[
\lb{RelaF2}
\begin{aligned}
\m_{\a,n}(p)=\n_{\a,n}(v), \qq  \n_{\a,n}(p)=\m_{\a,n}(v), \qq
\t_{\a,n}(p)=\vr_{\a,n}(v), \qq  \vr_{\a,n}(p)=\t_{\a,n}(v).
\end{aligned}
\]
These identities prove the relation in \er{RelaE} associated with $\1_-$.
Combining \er{RelaF} and \er{RelaF2}, we further obtain the following relations:
$$
\begin{aligned}
\gr_{\a,n}(p)=-\ln|(s_\a \vp_1+c_\a \vp_2)(1,\m_{\a,n}(p),p)|
=-\ln|(c_\a \vt_1-s_\a \vt_2)(1,\n_{\a,n}(v),v)|=\gs_{\a,n}(v),
\\
\gs_{\a,n}(p)=-\ln|(c_\a \vt_1-s_\a \vt_2)(1,\n_{\a,n}(p),p)|
=-\ln|(s_\a \vp_1+c_\a \vp_2)(1,\m_{\a,n}(v),v)|=\gr_{\a,n}(v),
\\
\gt_{\a,n}(p)=-\ln|(c_\a \vp_1-s_\a \vp_2)(1,\t_{\a,n}(p),p)|
=-\ln|(s_\a\vt_1+c_\a\vt_2)(1,\vr_{\a,n}(v),v)|=\gu_{\a,n}(v),
\\
\gu_{\a,n}(p)=-\ln|(s_\a \vt_1+c_\a\vt_2)(1,\vr_{\a,n}(p),p)|
=-\ln|(c_\a \vp_1-s_\a\vp_2)(1,\t_{\a,n}(v),v)|=\gt_{\a,n}(v).
\end{aligned}
$$
Since these relations hold for every $n\in\Z$, the first identity in
\er{RelaN} follows immediately.

We next establish the identities associated with $\mU_1$. Define
$$
\mU_1=\ma c_{2\a} & -s_{2\a} \\ -s_{2\a} & -c_{2\a} \am\cR,  \qqq
U_{\a}=\ma c_\a & -s_\a \\ -s_\a & -c_\a \am\cR.
$$
The gauge transformation of the system \er{eqZS} yields
\[
\lb{eqU1}
\tes U_\a(J\frac{d}{dx}+V)U_\a^*f=J f'+W f=\l f,
\]
where
$$
W=\ma w_1 & w_2 \\ w_2 & -w_1\am=U_\a V U_\a^*,\qqq \text{and}  \qqq
w=\mU_1 v.
$$
Let $Y$ denote the normalized fundamental matrix solution of
\er{eqU1}. By standard results in ODE theory, the fundamental
matrices corresponding to \(v\) and \(w\) are related by
\[
\lb{RelaF1}
\begin{aligned}
&Y(1,\l,w) =\ma -c_\a & s_\a \\ s_\a & c_\a\am Y^{-1}\ma -c_\a &
s_\a \\ s_\a & c_\a\am (1,\l,v).
\end{aligned}
\]
Taking suitable linear combinations of the entries in \er{RelaF1},
we obtain the following identities:
$$
\begin{aligned}
(c_\a\vp_1-s_\a\vp_2)(1,\l,w)=(c_\a\vp_1-s_\a\vp_2)(1,\l,v),
\\
(s_\a\vt_1+c_\a\vt_2)(1,\l,w)=(s_\a\vt_1+c_\a\vt_2)(1,\l,v),
\\
(s_\a\vp_1+c_\a\vp_2)(1,\l,w)=(c_\a\vt_1-s_\a\vt_2)(1,\l,v),
\\
(c_\a\vt_1-s_\a\vt_2)(1,\l,w)=(s_\a\vp_1+c_\a\vp_2)(1,\l,v),
\end{aligned}
$$
which implies that the associated eigenvalues satisfy
\[
\lb{RelaE1}
\begin{aligned}
(\m_{\a,n},\n_{\a,n},\t_{\a,n},\vr_{\a,n})\circ\mU_1
=(\m_{\a,n},\n_{\a,n},\vr_{\a,n},\t_{\a,n}), \qqq \forall \ n\in\Z.
\end{aligned}
\]
Recall that the Wronskian satisfies $\vt_1\vp_2-\vt_2\vp_1=1$.
Combining this relation with \er{RelaF1}, we obtain, for every \(n\in\Z\),
\[
\lb{RelaN1}
\begin{aligned}
e^{-\gr_{\a,n}(w)} &=|\big(s_\a\vp_1+c_\a\vp_2\big)(1,\m_{\a,n}(w),w)|
=|\big(c_\a\vt_1-s_\a \vt_2\big)(1,\m_{\a,n}(v),v)|
\\
&=|\big(s_\a\vp_1+c_\a \vp_2\big)(1,\m_{\a,n}(v),v)|^{-1}
=e^{\gr_{\a,n}(v)}.
\end{aligned}
\]
Similarly, we obtain
\[
\lb{RelaN2}
\begin{aligned}
e^{-\gs_{\a,n}(w)}
=|\big(c_\a\vt_1-s_\a\vt_2\big)(1,\n_{\a,n}(w),w)|
=|\big(c_\a\vt_1-s_\a\vt_2\big)(1,\n_{\a,n}(v),v)|^{-1}
=e^{\gs_{\a,n}(v)},
\\
e^{-\gt_{\a,n}(w)}
=|\big(c_\a\vp_1-s_\a\vp_2\big)(1,\t_{\a,n}(w),w)|
=|\big(s_\a\vt_1+c_\a  \vt_2\big)(1,\vr_{\a,n}(v),v)|^{-1}
=e^{\gu_{\a,n}(v)},
\\
e^{-\gu_{\a,n}(w)}
=|\big(s_\a\vt_1+c_\a\vt_2\big)(1,\vr_{\a,n}(w),w)|
=|\big(c_\a\vt_1-s_\a\vt_2\big)(1,\t_{\a,n}(v),v)|^{-1}
=e^{\gt_{\a,n}(v)}.
\end{aligned}
\]
Combining \er{RelaE1}--\er{RelaN2}, we establish all the identities
for $\mU_1$. Furthermore, \er{RelaF2} and \er{RelaE1} imply that,
for all $n\in\Z$,
$$
\begin{aligned}
\m_{\a,n}\ci(\1_-\mU_1)=\n_{\a,n}\ci\mU_1=\n_{\a,n}, \qqq
\n_{\a,n}\ci(\1_-\mU_1)=\m_{\a,n}\ci\mU_1=\m_{\a,n},
\\
\t_{\a,n}\ci(\1_-\mU_1)=\vr_{\a,n}\ci\mU_1=\t_{\a,n}, \qqq
\vr_{\a,n}\ci(\1_-\mU_1)=\t_{\a,n}\ci\mU_1=\vr_{\a,n},
\end{aligned}
$$
and
$$
\begin{aligned}
\gr_{\a,n}\ci(\1_-\mU_1)=\gs_{\a,n}\ci\mU_1=-\gs_{\a,n}, \qqq
\gs_{\a,n}\ci(\1_-\mU_1)=\gr_{\a,n}\ci\mU_1=-\gr_{\a,n},
\\
\gt_{\a,n}\ci(\1_-\mU_1)=\gu_{\a,n}\ci\mU_1=-\gt_{\a,n}, \qqq
\gu_{\a,n}\ci(\1_-\mU_1)=\gt_{\a,n}\ci\mU_1=-\gu_{\a,n}.
\end{aligned}
$$
Since $\mU_2=\1_-\mU_1$, these identities establish \er{RelaE} and
\er{RelaN} for $\mU_2$, completing the proof.           \BBox

\begin{corollary} \lb{Ce+n}
Let $\a\in\A$ and $F_{2\a}, \mU_1, \mU_2$ be defined by \er{DefF_a},
\er{DefF_12}.

\no i) All two-spectra mappings $\t_\a\star \m_\a$,
$\t_\a\star \n_\a $, $\vr_\a\star \m_\a$, $\vr_\a\star \n_\a $
are isomorphic and satisfy
\[
\lb{2sma'}
\begin{aligned}
\cS_\a(\t\star \m)\ci F_{2\a}^*=\t_\a\star \m_\a
=(\vr_\a\star \n_\a )\circ\1_- \\
=(\vr_\a\star \m_\a)\circ\mU_1=(\t_\a\star \n_\a )\circ\mU_2.
\end{aligned}
\]

\no ii) All mappings $\cS_\a\m\ts \gr$, $\m_\a\ts (\pm\gr_\a )$,
$\n_\a  \ts (\pm\gs_\a )$ are isomorphic. Likewise, the mappings
$\cS_\a\t\ts \gt$, $\vr_\a\ts(\pm\gu_\a )$ and $\t_\a\ts(\pm\gt_\a
)$ are isomorphic. Moreover, they satisfy
\[
\lb{SpN1'}
\begin{aligned}
(\cS_\a\m\ts \gr) \ci F_{2\a}^*&=\m_\a\ts \gr_\a
=(\n_\a \ts \gs_\a  )\circ \1_-
\\
&=(\m_\a \ts (-\gr_\a ))\circ \mU_1
=(\n_\a  \ts (-\gs_\a ))\circ \mU_2,
\end{aligned}
\]
and
\[
\lb{SpN2'}
\begin{aligned}
(\cS_\a\t\ts \gt) \ci F_{2\a}^*&= \t_\a\ts \gt_\a
=(\vr_\a \ts \gu_\a  )\circ \1_-
\\
&=(\vr_\a \ts (-\gu_\a ) )\circ \mU_1
=(\t_\a \ts (-\gt_\a ))\circ \mU_2.
\end{aligned}
\]
\end{corollary}

\no {\bf{Proof.}} i) Lemmas \ref{LrotationF_a} and
\ref{LrotationF_12}  imply that the two-spectra mappings
$\cS_\a(\t\star\m)$, $\t_\a\star \m_\a$, $\t_\a\star \n_\a $,
$\vr_\a\star \m_\a$, $\vr_\a\star \n_\a $ are isomorphic and satisfy
\er{2sma'}.

ii) The same lemmas show that the mappings $\cS_\a\m\ts \gr,\m_\a\ts
(\pm\gr_\a )$, $\n_\a\ts (\pm\gs_\a )$ are pairwise isomorphic, as
are the mappings $\cS_\a\t\ts \gt$, $\vr_\a\ts(\pm\gu_\a )$ and
$\t_\a\ts(\pm\gt_\a )$. Furthermore, all these mappings satisfy
\er{SpN1'} and \er{SpN2'}.                \BBox

\section{Isomorphic inverse problems on the unit interval \lb{Sec4}}
\setcounter{equation}{0}

\subsection{Proof of main theorems}
We now transfer the basic results through the gauge transformation.
The same spectral translation and the same unitary change of
potential will be used in each argument.

\no {\bf Proof of Theorem \ref{T4sma}.} Let $v\in\mH$. Recall  that
$F_{2\a}$ and $\cS_\a$ are defined in \er{DefF_a} and
\er{DefShiftO}, respectively. From \er{cA} in Lemma
\ref{LrotationF_a}, we obtain
$$
\begin{aligned}
\tes \gf_{\a,2n-1}\ci F_{2\a}=\frac{1}{2}(\vr_{\a,n}-\t_{\a,n})
\ci F_{2\a}=\frac{1}{2}\big((\vr_n-\a)-(\t_n-\a)\big)
=\frac{1}{2}(\vr_n-\t_n),
\\
\tes \gf_{\a,2n}\ci F_{2\a}=\frac{1}{2}(\n_{\a,n}-\m_{\a,n})
\ci F_{2\a}=\frac{1}{2}\big((\n_n-\a)-(\m_n-\a)\big)
=\frac{1}{2}(\n_n-\m_n).
\end{aligned}
$$
These identities yield $\gf_\a \ci F_{2\a}=\gf$. Since the mapping
$\gf: \mH\to \ell^2$ is a RAB, Theorem \ref{L4sm}  and the preceding
relation imply that the mapping $\gf_\a : \mH\to \ell^2$ is also a
RAB. Using the norm-preserving property of \(F_{2\alpha}\) and
applying \er{egf}, we obtain the following estimates:
$$
\begin{aligned}
&\|\gf_\a (u)\|=\|\gf(v)\|\le\sqrt2 \|v\|=\sqrt2\|u\|,
\\
\|u\|=\|v\|&\le 2\|\gf(v)\|\big(1+\|\gf(v)\|\big)=2\|\gf_\a (u)
\|\big(1+\|\gf_\a (u)\|\big),
\end{aligned}
$$
where $u=F_{2\a} v$. These establish \er{egf1}.            \BBox

\begin{theorem}
\label{TSsa} Let $\a\in \A$.

\no i) The mapping $(\cS_\a\m)\ts \gr : \mH\to \gJ^o_\a\ts \ell^2$
is a RAB from $\mH$ onto $\gJ^o_\a\ts \ell^2$.

\no ii) The mapping $(\cS_\a\t) \ts \gt: \mH \to \gJ^1_\a\ts \ell^2$
is a RAB from $\mH$ onto $\gJ^1_\a\ts \ell^2$.

\no iii) The mapping $\cS_\a\m:\mH_{eo}\to\gJ^o_\a$ is a RAB
from $\mH_{eo}$ onto $\gJ^o_\a$.
\end{theorem}

\no {\bf Proof.} In standard sequence coordinates, $\cS_\a$ is an
analytic isomorphism. Parts i) and ii) therefore follow by composing
the RABs in Theorem \ref{Liso} with $\cS_\a\times I$, where $I$ is
the identity on $\ell^2$. Part iii) follows by composing the RAB in
Theorem \ref{L1seo} with $\cS_\a$. \BBox

\no {\bf Proof of Theorem \ref{Tiso2sm}.} By Theorem \ref{Liso}, the
two-spectra mapping $\t\star\m$ is a RAB from $\mH$ onto $\gJ$.
Since $\cS_\a$ is a real-analytic isomorphism from $\gJ$ onto
$\gJ_\a$, the composition $\cS_\a( \t\star\m)$ is also a RAB from
$\mH$ onto $\gJ_\a$. Moreover, Corollary \ref{Ce+n} shows that all
two-spectra mappings $\cS_\a( \t\star\m)$, $\t_\a\star\m_\a$,
$\t_\a\star\n_\a $, $\vr_\a\star\m_\a$, $\vr_\a\star\n_\a $, which
map $\mH$ to $\gJ_\a$, are isomorphic. Consequently, each of these
mappings is a RAB from $\mH$ onto $\gJ_\a$ and satisfies \er{2sma}.
\BBox

\no {\bf Proof of Theorem \ref{Tisoe+n}.} i) From Theorem
\ref{TSsa}, the mapping  $(\cS_\a\m)\ts\gr: \mH\to
\gJ^o_\a\ts\ell^2$ is a RAB from $\mH$ onto $\gJ^o_\a\ts\ell^2$.
Corollary \ref{Ce+n} further shows that the mappings
$(\cS_\a\m)\ts\gr$, $\m_\a\ts(\pm\gr_\a )$, $\n_\a \ts(\pm\gs_\a )$
from $\mH$ to $\gJ^o_\a\ts\ell^2$ are isomorphic. Consequently, each
of these mappings is a RAB from $\mH$ onto $\gJ^o_\a\ts\ell^2$ and
satisfies \er{SpN1}.

ii) Similarly, Theorem \ref{TSsa} implies that the mapping
$(\cS_\a\t)\ts\gt: \mH\to \gJ^1_\a\ts\ell^2$ is a RAB from $\mH$
onto $\gJ^1_\a\ts\ell^2$. By Corollary \ref{Ce+n}, all the mappings
$(\cS_\a\t)\ts\gt$, $\t_\a\ts(\pm\gt_\a )$, $\vr_\a\ts(\pm\gu_\a )$
from $\mH$ to $\gJ^1_\a\ts\ell^2$ are isomorphic. Thus, each of
these mappings is a RAB from $\mH$ onto $\gJ^1_\a\ts\ell^2$ and
satisfies \er{SpN2}.                 \BBox

\no {\bf Proof of Corollary \ref{Ceoa}.}
Recall that the spaces $\mH_{eo}$ and $\mH_\a$ are defined in \er{def-ha}.
We show that $F_{2\a}:\mH_{eo}\to\mH_\a$ is a bijection. For any $v\in\mH_{eo}$ we have
$$
u(x)=F_{2\a} v(x)=\ma c_{2\a x}v_1+s_{2\a x}v_2
\\
-s_{2\a x}v_1+c_{2\a x}v_2 \am,
$$
By the definition of $\mH_{eo}$, we obtain
$$
u(1-x)=\ma c_{2\a} u_1-s_{2\a} u_2 \\ -s_{2\a} u_1-c_{2\a} u_2 \am(x)=e^{2\a J}J_1u(x).
$$
Hence $u\in\mH_\a$. A direct computation shows that, for every
$u\in\mH_\a$, there exists a unique $v\in\mH_{eo}$ such that
$u=F_{2\a} v$.

From Theorems \ref{Liso} and \ref{L1seo}, we deduce that the one--spectra mappings
$$
\cS_\a\m, \cS_\a\n: \mH_{eo}\to \gJ^o_\a, \qqq \cS_\a\t, \cS_\a\vr: \mH_{eo}\to\gJ^1_\a,
$$
are all RABs. Combining these facts with \er{cA} yields identity \er{shifteo}. \BBox

\no {\bf Proof of Corollary \ref{Creplacing}.} Recall that the
mapping  $\cU_{\a,\s}$, defined in \er{defUa}, is a RAB from $\mH$
onto itself. By Theorems \ref{Llamp}, \ref{Lreplacing}, and
Lemma \ref{LrotationF_a}, we deduce that
$$
\big(\z_\a \star\phi_\a \big)\circ\cU_{\a,\s}= \big(\z_\a
\star\phi_\a \big)\circ(F_{2\a}\ci \cU_\s\ci
F_{2\a}^*)=\t_\a\star\m_\a,
$$
which yields \er{rma}. By Theorem \ref{Tiso2sm}, the  two-spectra
mapping $\t_\a\star\m_\a$ is a RAB from $\mH$ onto $\gJ_\a$. Thus
the mapping $\z_\a \star\phi_\a $ is also a RAB from $\mH$ onto
$\gJ_\a$.             \BBox

\no {\bf Proof of Theorem \ref{Tcanonicala}.}
Let $a, b: \mH\to\R$ be continuously Fr\'echet differentiable functions.
Let $F:\mH\to\mH$ be a linear multiplication operator of the form
$$
Fv=F(x)v(x),
$$
where $F(x)$ is a $2\ts 2$ matrix.
Define the composite functions $A(v)=a(Fv)$ and $B(v)=b(Fv)$ for  $v\in\mH$.
Writing gradients as row vectors and setting $u=F v$, the chain rule gives
\[
\lb{FSz}
\begin{aligned}
(A_v',B_v')_0= (A_{v_1}',A_{v_2}')^\top J (B_{v_1}',B_{v_2}')=
(a_{u_1}',a_{u_2}')^\top F^\top J F (b_{u_1}',b_{u_2}').
\end{aligned}
\]
In particular, if \(FJF^{\top}=J\), then
$$
A(v)\wedge B(v)=a(u)\wedge b(u).
$$
By Lemma \ref{LrotationF_a}, we deduce that $\cS_\a\m=\m_\a \circ
F_{2\a}$. Combining this identity with Theorem \ref{Lcanonical},
identity \er{FSz},  and the fact that $F_{2\a}^* J F_{2\a}=J$, we
obtain, for all $n,j\in\Z$,
$$
\begin{aligned}
 \m_{\a,n}(u)\we \m_{\a,j}(u)=\m_n(v)\we \m_j(v)=0,  \qqq
 \gr_{\a,n}(u)\we \gr_{\a,j}(u)=\gr_n(v)\we \gr_j(v)=0,  \\
 \gr_{\a,n}(u)\we \m_{\a,j}(u)=\gr_n(v)\we \m_j(v)=\d_{n,j},
 \end{aligned}
$$
which gives the first line in \er{cv1a}. The remaining identities of
\er{cv1a} follow by analogous arguments;  see also \cite{KZ26}. By a
well-known result in functional analysis (see \cite{PT87}), the
gradient family $\{ \m_{\a,n}', \gr'_{\a,n}\}_{n\in \Z}$, forms a
basis of $\mH$, as do the gradient families associated with the
other spectral mappings. This completes the proof.  \BBox

\subsection{Normalizing mappings }
We introduce the normalizing mappings $\ga_\a =(\ga_{\a,n})_{n\in\Z}$,
$\gb_\a =(\gb_{\a,n})_{n\in\Z}$, $\gc_\a =(\gc_{\a,n})_{n\in\Z}$,
$\gd_\a =(\gd_{\a,n})_{n\in\Z}$, with components defined by
\[
\lb{DefNormalizing}
\begin{aligned}
e^{-\ga_{\a,n}}=\|\vp(\cdot,\m_{\a,n},v)\|^2,\qqq
e^{-\gb_{\a,n}}=\|\vt(\cdot,\n_{\a,n},v)\|^2,\\
e^{-\gc_{\a,n}}=\|\vp(\cdot,\t_{\a,n},v)\|^2,\qqq
e^{-\gd_{\a,n}}=\|\vt(\cdot,\vr_{\a,n},v)\|^2.
\end{aligned}\]

\begin{lemma}
\lb{Lrotanla} Let $\a\in\A$ and $F_{2\a}$ be defined by \er{DefF_a}.
Then  the normalizing mappings $\ga_\a $, $\gb_\a $, $\gc_\a $,
$\gd_\a $, defined by \er{DefNormalizing}, satisfy
\[
\lb{estin}
(\ga_\a,\gb_\a,\gc_\a,\gd_\a)\ci F_{2\a}=(\ga,\gb,\gc,\gd).
\]
\end{lemma}

\no {\bf Proof.} Let $u=F_{2\a}v$ for $v\in\mH$.
Lemma~\ref{LrotationF_a}  gives, for example,
$\vp(x,\m_{\a,n}(u),u)=F_\a(x)\vp(x,\m_n(v),v)$. These
eigenfunctions have the same norm, i.e.
$$
\begin{aligned}
\ga_{\a,n}(u)=-\ln\int^1_0\vp^2(x,\m_{\a,n},u)dx= -\ln\int^1_0
\big(F_\alpha(x)\varphi(x,\mu_n,v)\big)^2dx=\ga_n(v), \qq \forall~
n\in\Z.
\end{aligned}
$$
By the same argument with the other three eigenfunctions,  we obtain
the following identities valid for all $n\in\Z$,
$$
\begin{aligned}
\gb_{\a,n}(u)=-\ln\int^1_0\vt^2(x,\n_{\a,n},u)dx=
-\ln\int^1_0 \big(F_\a(x)\vt(x,\n_n,v)\big)^2dx =\gb_n(v),
\\
\gc_{\a,n}(u)=-\ln\int^1_0\vp^2(x,\t_{\a,n},u)dx=
-\ln\int^1_0 \big(F_\a(x)\vp(x,\t_n,v)\big)^2dx=\gc_n(v),
\\
\gd_{\a,n}(u)=-\ln\int^1_0\vt^2(x,\vr_{\a,n},u)dx=
-\ln\int^1_0 \big(F_\a(x)\vt(x,\vr_n,v)\big)^2dx=\gd_n(v).
\end{aligned}
$$
These prove the identity \er{estin}. \BBox

We express the normalizing constants $\ga_n,\gb_n $,
 $\gc_n, \gd_n$ (for $n\in\Z$) in terms of the norming
  constants $\gr_n, \gs_n $, $\gt_n, \gu_n$ in the following form:
\begin{equation}
\label{ss1}
\begin{aligned}
\tes
& e^{-\ga_n }=-\dot \vp _1(1, \m_n)\vp_2(1,\m_n)
=e^{ \ga_n ^\bu -\gr_n },\qqq  \ga_n ^\bu=\ln |\dot \vp _1(1, \m_n)|, \\
& e^{-\gb_n}=\dot \vt_2(1,\n_n)\vt_1(1,\n_n)
=e^{\gb_n^\bu -\gs_n },\qqq \gb_n^\bu=\ln |\dot \vt_2(1,\n_n)|, \\
& e^{-\gc_n}=\dot \vp_2(1, \t_n)\vp_1(1,\t_n)=e^{\gc_n^\bu-\gt_n },
\qqq \gc_n^\bu=\ln |\dot \vp_2(1, \t_n)|, \\
& e^{-\gd_n}=-\dot \vt_1(1,\vr_n)\vt_2(1,\vr_n)
=e^{  \gd_n^\bu   -\gu_n },\qqq \gd_n^\bu=\ln |\dot \vt_1(1,\vr_n)|,
\end{aligned}
\end{equation}
where $\dot u=\frac{\pa u}{\l}$.
We establish the relations between the norming and normalizing mappings.

\begin{lemma} \lb{NandN}
Let $\a\in\A$. Then for any $n\in\Z$, the normalizing constants
defined by  \er{DefNormalizing} satisfy
\[
\lb{nlc}
\begin{aligned}
&e^{-\ga_{\a,n}}= \int^1_0|\vp(x,\m_{\a,n})|^2dx=e^{A_n^\bu-\gr_{\a,n}}, \qqq
A_n^\bu=\ln|c_\a \dot\vp_1(1,\m_{\a,n})-
s_\a\dot\vp_2(1,\m_{\a,n})|, \\
&e^{-\gb_{\a,n}}= \int^1_0|\vt(x,\n_{\a,n})|^2dx=e^{B_n^\bu-\gs_{\a,n}}, \qqq
B_n^\bu=\ln|s_\a\dot\vt_1(1,\n_{\a,n})+c_\a \dot
\vt_2(1,\n_{\a,n}) |,\\
&e^{-\gc_{\a,n}}= \int^1_0|\vp(x,\t_{\a,n})|^2dx=e^{C_n^\bu-\gt_{\a,n}}, \qqq
C_n^\bu=\ln|s_\a\dot\vp_1(1,\t_{\a,n})+c_\a \dot
\vp_2(1,\t_{\a,n})|,\\
&e^{-\gd_{\a,n}}=\int^1_0|\vt(x,\vr_{\a,n})|^2dx
=e^{D_n^\bu-\gu_{\a,n}}, \qqq
D_n^\bu=\ln|c_\a \dot\vt_1(1,\vr_{\a,n})-s_\a\dot
\vt_2(1,\vr_{\a,n})|.
\end{aligned}
\]
\end{lemma}

\no {\bf Proof.} Let $u=F_{2\a} v$ for $v\in\mH$. It follows from \er{3.1} that
$$
\big(c_\a\vp_1-s_\a\vp_2\big)(1,\z,u)=\vp_1(1,\l,v), \qqq \z=\l-\a.
$$
Differentiating this identity with respect to the spectral parameter and using \er{cA}, we obtain
\[
\lb{driv1}
\big(c_\a \dot\vp_1-s_\a \dot\vp_2\big)(1,\m_{\a,n},u)=\dot\vp_1(1,\m_n,v), \qqq
\frac{\vp_2(1,\m_{\a,n},u)}{c_\a}=\vp_2(1,\m_n,v),\qqq \forall \ n\in\Z.
\]
Applying the same argument to the remaining entries of the
fundamental matrix  yields the following identities valid for all
$n\in\Z$,
\[
\lb{driv}
\begin{aligned}
&\big(s_\a \dot\vp_1+c_\a \dot\vp_2\big)(1,\t_{\a,n},u)=\dot\vp_2(1,\t_n,v),\qqq
\frac{\vp_1(1,\t_{\a,n},u)}{c_\a}=\vp_1(1,\t_n,v), \\
&\big(c_\a \dot\vt_1-s_\a \dot\vt_2\big)(1,\vr_{\a,n},u)=\dot\vt_1(1,\vr_n,v), \qqq
\frac{\vt_2(1,\vr_{\a,n},u)}{c_\a}=\vt_2(1,\vr_n,v), \\
&\big(s_\a \dot\vt_1+c_\a \dot\vt_2\big)(1,\n_{\a,n},u)=\dot\vt_2(1,\n_n,v), \qqq
\frac{\vt_1(1,\n_{\a,n},u)}{c_\a}=\vt_1(1,\n_n,v).
\end{aligned}
\]
Combining \er{ss1} with \er{estin}, \er{driv1}, and \er{driv}, we obtain
$$
\begin{aligned}
\ga_{\a, n}(u)=\ga_n(v)=& -\ln\Big| \frac{\vp_2(1,\m_{\a,n},u)}{c_\a}
\big(c_\a \dot\vp_1-s_\a \dot\vp_2\big)(1,\m_{\a,n},u)\Big| \\
=&-\ln\Big|\frac{\vp_2(1,\m_{\a,n},u)}{c_\a}\Big|
-\ln\Big|(c_\a \dot\vp_1-s_\a \dot\vp_2)(1,\m_{\a,n},u)\Big|,
\end{aligned}
$$
$$
\begin{aligned}
\gb_{\a, n}(u)=\gb_n(v)=& -\ln\Big|\frac{\vt_1(1,\n_{\a,n},u)}{c_\a}
\big(s_\a \dot\vt_1+c_\a \dot\vt_2\big)(1,\n_{\a,n},u)\Big| \\
=&-\ln\Big|\frac{\vt_1(1,\n_{\a,n},u)}{c_\a}\Big|
-\ln\Big|(s_\a \dot\vt_1+c_\a \dot\vt_2)(1,\n_{\a,n},u)\Big|, \\
\end{aligned}
$$
$$
\begin{aligned}
\gc_{\a, n}(u)=\gc_n(v)=& -\ln\Big| \frac{\vp_1(1,\t_{\a,n},u)}{c_\a}
\big(s_\a \dot\vp_1+c_\a \dot\vp_2\big)(1,\t_{\a,n},u)\Big| \\
=&-\ln\Big|\frac{\vp_1(1,\t_{\a,n},u)}{c_\a}\Big|
-\ln\Big|(s_\a \dot\vp_1+c_\a \dot\vp_2)(1,\t_{\a,n},u)\Big|, \\
\end{aligned}
$$
$$
\begin{aligned}
\gd_{\a, n}(u)=\gd_n(v)=& -\ln\Big| \frac{\vt_2(1,\vr_{\a,n},u)}{c_\a}
\big(c_\a \dot\vt_1-s_\a \dot\vt_2\big)(1,\vr_{\a,n},u)\Big| \\
=&-\ln\Big|\frac{\vt_2(1,\vr_{\a,n},u)}{c_\a}\Big|
-\ln\Big|(c_\a \dot\vt_1-s_\a \dot\vt_2)(1,\vr_{\a,n},u)\Big|.
\end{aligned}
$$
These identities prove \er{nlc}.              \BBox

\begin{lemma} Let $\a\in\A$.
The normalizing mappings $\ga_\a, \gb_\a$, $\gc_\a$, $\gd_\a$,
defined in \er{DefNormalizing}, satisfy
\[
\lb{nlca}
(\ga_\a ,\gb_\a ,\gc_\a ,\gd_\a )
=(\gb_\a ,\ga_\a ,\gd_\a ,\gc_\a )\ci\1_-.
\]
\end{lemma}

\no {\bf Proof.} Using \er{nlc}, we deduce that, for all $n\in\Z$,
the normalizing constants satisfy
$$
\ga_{\a,n}=\gr_{\a,n}-A_n^\bu, \qq
\gb_{\a,n}=\gs_{\a,n}-B_n^\bu, \qq
\gc_{\a,n}=\gt_{\a,n}-C_n^\bu, \qq
\gd_{\a,n}=\gu_{\a,n}-D_n^\bu,
$$
Let $p=-v$. Differentiating the identities in \er{RelaF} with
respect  to the spectral parameter, we obtain
\[
\lb{buA}
\begin{aligned}
(c_\a  \dot\vp_1-s_\a  \dot\vp_2)(1,\l,p)
=-(c_\a  \dot\vt_2+s_\a  \dot\vt_1)(1,\l,v), \\
(s_\a  \dot\vt_1+c_\a  \dot\vt_2)(1,\l,p)
=-(c_\a  \dot\vp_1-s_\a  \dot\vp_2)(1,\l,v), \\
(s_\a  \dot\vp_1+c_\a  \dot\vp_2)(1,\l,p)
=(c_\a  \dot\vt_1-s_\a  \dot\vt_2)(1,\l,v), \\
(c_\a  \dot\vt_1-s_\a  \dot\vt_2)(1,\l,p)
=(s_\a  \dot\vp_1+c_\a  \dot\vp_2)(1,\l,v).
\end{aligned}
\]
Combining \er{buA} with \er{RelaE}, we obtain the following identities for all $n\in\Z$:
$$
\begin{aligned}
A_n^\bu(p)=\ln|\big(c_\a\dot\vp_1-s_\a\dot\vp_2\big) (1,\m_{\a,n}(p),p)|
=\ln|\big(s_\a\dot\vt_1+c_\a\dot\vt_2\big)(1,\n_{\a,n}(v),v)| =B_n^\bu(v),
\end{aligned}
$$
$$
\begin{aligned}
B_n^\bu(p)=\ln|\big(s_\a\dot\vt_1
+c_\a\dot\vt_2\big)(1,\n_{\a,n}(p),p)|
=\ln|\big(c_\a\dot\vp_1-s_\a\dot\vp_2\big) (1,\m_{\a,n}(v),v)|=A_n^\bu(v),
\end{aligned}
$$
$$
\begin{aligned}
C_n^\bu(p)=\ln|\big(s_\a\dot\vp_1+c_\a\dot\vp_2\big) (1,\t_{\a,n}(p),p)|
=\ln|\big(c_\a\dot\vt_1
-s_\a\dot\vt_2\big)(1,\vr_{\a,n}(v),v)|=D_n^\bu(v),
\end{aligned}
$$
$$
\begin{aligned}
D_n^\bu(p)=\ln|\big(c_\a\dot\vt_1
-s_\a\dot\vt_2\big)(1,\vr_{\a,n}(p),p)|
=\ln|\big(s_\a  \dot\vp_1
+c_\a\dot\vp_2\big)(1,\t_{\a,n}(v),v)|=C_n^\bu(v).
\end{aligned}
$$
Combining this observation with \er{RelaN} and the spectral  swaps
in \er{RelaF2}, we conclude that
$$
\begin{aligned}
&\ga_{\a,n}\circ\1_-=\gs_{\a,n}-B_n^\bu=\gb_{\a,n}, \qqq
\gb_{\a,n}\circ\1_-=\gr_{\a,n}-A_n^\bu=\ga_{\a,n}, \qq \forall n\in\Z, \\
&\gc_{\a,n}\circ\1_-=\gu_{\a,n}-D_n^\bu=\gd_{\a,n}, \qqq
\gd_{\a,n}\circ\1_-=\gt_{\a,n}-C_n^\bu=\gc_{\a,n}, \qq \forall n\in\Z,
\end{aligned}
$$
These identities proves \er{nlca}.     \BBox

\begin{theorem}
\label{Te+nla} Let $\a\in\A$. Then

\no i) The mappings $(\cS_\a\m)\ts \ga$, $\m_\a\ts \ga_\a $ and
$\n_\a \ts\gb_\a $ are bijections from $\mH$ onto
$\gJ^o_\a\ts\ell^2$.  Moreover, they satisfy
\[
\lb{eqe+nl-1}
\m_\a\ts\ga_\a =(\n_\a \ts\gb_\a )\circ\1_-=(\cS_\a\m\ts\ga)\circ F^*_{2\a}.
\]
\no ii) Each of the mappings $(\cS_\a\t)\ts\gc$, $\t_\a\ts\gc_\a$
and  $\vr_\a\ts\gd_\a $ is a bijection from $\mH$ onto $\gJ^1_\a\ts
\ell^2$. Moreover, they satisfy
$$
\t_\a\ts\gc_\a =(\vr_\a\ts\gd_\a )\circ\1_-=(\cS_\a\t\ts\gc)\circ F^*_{2\a}.
$$
\end{theorem}

\no {\bf Proof.} We first prove i). By Theorem \ref{Le+nl},
\(\m\ts\ga\)  is a bijection from \(\mH\) onto
\(\gJ^o\times\ell^2\). Since $\cS_\a$ is bijective, the mapping
$(\cS_\a\m)\times\ga$ is therefore a bijection from $\mH$ to
$\gJ^o_\a\ts\ell^2$. Combining the spectral identities \er{cA},
\er{RelaE} with the normalizing identities \er{estin} and \er{nlca},
we obtain \er{eqe+nl-1}. Furthermore, Lemmas \ref{LrotationF_a},
\ref{LrotationF_12}, and \ref{Lrotanla} imply that the mappings
$\cS_\a\m\ts \ga$, $\m_\a\ts \ga_\a $ and $\n_\a \ts\gb_\a $ are
bijections from $\mH$ onto $\gJ^o_\a\ts\ell^2$.

Part ii) follows by an analogous argument.           \BBox

\no \textbf{Acknowledgments.}  {\footnotesize The research of E.
Korotyaev is supported  by PJSC "Gazprom Neft",  aggr.
GPN-26/03000/00502/R. Various parts of this paper were written at
Renmin University of China, Beijing. E. K. thanks the university for
its hospitality.}

\end{document}